\documentclass{article}
\usepackage{amssymb,amsmath,theorem,euscript}

\input{epsf}

\newcounter{sec}

\def\sm{\smallskip}

\newcounter{punct}[sec]

\def\punct{\refstepcounter{punct}{\arabic{sec}.\arabic{punct}.  }}

\def\COUNTERS{\addtocounter{sec}{1}
              \setcounter{punct}{0}
          \setcounter{equation}{0}
          \setcounter{theorem}{0}
                  }

\newtheorem{theorem}{Theorem}[sec]
\newtheorem{proposition}[theorem]{Proposition}
\newtheorem{lemma}[theorem]{Lemma}
\newtheorem{corollary}[theorem]{Corollary}
\newtheorem{observation}[theorem]{Observation}

\begin{document}

 \def\ov{\overline}
\def\wt{\widetilde}
 \newcommand{\rk}{\mathop {\mathrm {rk}}\nolimits}
\newcommand{\Aut}{\mathop {\mathrm {Aut}}\nolimits}
\newcommand{\Out}{\mathop {\mathrm {Out}}\nolimits}
 \newcommand{\tr}{\mathop {\mathrm {tr}}\nolimits}
  \newcommand{\diag}{\mathop {\mathrm {diag}}\nolimits}
  \newcommand{\supp}{\mathop {\mathrm {supp}}\nolimits}
 
\renewcommand{\Re}{\mathop {\mathrm {Re}}\nolimits}

\def\Br{\mathrm {Br}}

\def\SL{\mathrm {SL}}
\def\SU{\mathrm {SU}}
\def\GL{\mathrm {GL}}
\def\U{\mathrm U}
\def\OO{\mathrm O}
 \def\Sp{\mathrm {Sp}}
 \def\SO{\mathrm {SO}}
\def\SOS{\mathrm {SO}^*}
 \def\Diff{\mathrm{Diff}}
 \def\Vect{\mathfrak{Vect}}
\def\PGL{\mathrm {PGL}}
\def\PU{\mathrm {PU}}
\def\PSL{\mathrm {PSL}}
\def\Symp{\mathrm{Symp}}
\def\End{\mathrm{End}}
\def\Mor{\mathrm{Mor}}
\def\Aut{\mathrm{Aut}}
 \def\PB{\mathrm{PB}}
 \def\cA{\mathcal A}
\def\cB{\mathcal B}
\def\cC{\mathcal C}
\def\cD{\mathcal D}
\def\cE{\mathcal E}
\def\cF{\mathcal F}
\def\cG{\mathcal G}
\def\cH{\mathcal H}
\def\cJ{\mathcal J}
\def\cI{\mathcal I}
\def\cK{\mathcal K}
 \def\cL{\mathcal L}
\def\cM{\mathcal M}
\def\cN{\mathcal N}
 \def\cO{\mathcal O}
\def\cP{\mathcal P}
\def\cQ{\mathcal Q}
\def\cR{\mathcal R}
\def\cS{\mathcal S}
\def\cT{\mathcal T}
\def\cU{\mathcal U}
\def\cV{\mathcal V}
 \def\cW{\mathcal W}
\def\cX{\mathcal X}
 \def\cY{\mathcal Y}
 \def\cZ{\mathcal Z}
\def\0{{\ov 0}}
 \def\1{{\ov 1}}
 \def\frA{\mathfrak A}
 \def\frB{\mathfrak B}
\def\frC{\mathfrak C}
\def\frD{\mathfrak D}
\def\frE{\mathfrak E}
\def\frF{\mathfrak F}
\def\frG{\mathfrak G}
\def\frH{\mathfrak H}
\def\frI{\mathfrak I}
 \def\frJ{\mathfrak J}
 \def\frK{\mathfrak K}
 \def\frL{\mathfrak L}
\def\frM{\mathfrak M}
 \def\frN{\mathfrak N} \def\frO{\mathfrak O} \def\frP{\mathfrak P} \def\frQ{\mathfrak Q} \def\frR{\mathfrak R}
 \def\frS{\mathfrak S} \def\frT{\mathfrak T} \def\frU{\mathfrak U} \def\frV{\mathfrak V} \def\frW{\mathfrak W}
 \def\frX{\mathfrak X} \def\frY{\mathfrak Y} \def\frZ{\mathfrak Z} \def\fra{\mathfrak a} \def\frb{\mathfrak b}
 \def\frc{\mathfrak c} \def\frd{\mathfrak d} \def\fre{\mathfrak e} \def\frf{\mathfrak f} \def\frg{\mathfrak g}
 \def\frh{\mathfrak h} \def\fri{\mathfrak i} \def\frj{\mathfrak j} \def\frk{\mathfrak k} \def\frl{\mathfrak l}
 \def\frm{\mathfrak m} \def\frn{\mathfrak n} \def\fro{\mathfrak o} \def\frp{\mathfrak p} \def\frq{\mathfrak q}
 \def\frr{\mathfrak r} \def\frs{\mathfrak s} \def\frt{\mathfrak t} \def\fru{\mathfrak u} \def\frv{\mathfrak v}
 \def\frw{\mathfrak w} \def\frx{\mathfrak x} \def\fry{\mathfrak y} \def\frz{\mathfrak z} \def\frsp{\mathfrak{sp}}
 \def\bfa{\mathbf a} \def\bfb{\mathbf b} \def\bfc{\mathbf c} \def\bfd{\mathbf d} \def\bfe{\mathbf e} \def\bff{\mathbf f}
 \def\bfg{\mathbf g} \def\bfh{\mathbf h} \def\bfi{\mathbf i} \def\bfj{\mathbf j} \def\bfk{\mathbf k} \def\bfl{\mathbf l}
 \def\bfm{\mathbf m} \def\bfn{\mathbf n} \def\bfo{\mathbf o} \def\bfp{\mathbf p} \def\bfq{\mathbf q} \def\bfr{\mathbf r}
 \def\bfs{\mathbf s} \def\bft{\mathbf t} \def\bfu{\mathbf u} \def\bfv{\mathbf v} \def\bfw{\mathbf w} \def\bfx{\mathbf x}
 \def\bfy{\mathbf y} \def\bfz{\mathbf z} \def\bfA{\mathbf A} \def\bfB{\mathbf B} \def\bfC{\mathbf C} \def\bfD{\mathbf D}
 \def\bfE{\mathbf E} \def\bfF{\mathbf F} \def\bfG{\mathbf G} \def\bfH{\mathbf H} \def\bfI{\mathbf I} \def\bfJ{\mathbf J}
 \def\bfK{\mathbf K} \def\bfL{\mathbf L} \def\bfM{\mathbf M} \def\bfN{\mathbf N} \def\bfO{\mathbf O} \def\bfP{\mathbf P}
 \def\bfQ{\mathbf Q} \def\bfR{\mathbf R} \def\bfS{\mathbf S} \def\bfT{\mathbf T} \def\bfU{\mathbf U} \def\bfV{\mathbf V}
 \def\bfW{\mathbf W} \def\bfX{\mathbf X} \def\bfY{\mathbf Y} \def\bfZ{\mathbf Z} \def\bfw{\mathbf w}
 \def\R {{\mathbb R }} \def\C {{\mathbb C }} \def\Z{{\mathbb Z}} \def\H{{\mathbb H}} \def\K{{\mathbb K}}
 \def\N{{\mathbb N}} \def\Q{{\mathbb Q}} \def\A{{\mathbb A}} \def\T{\mathbb T} \def\P{\mathbb P} \def\G{\mathbb G}
 \def\bbA{\mathbb A} \def\bbB{\mathbb B} \def\bbD{\mathbb D} \def\bbE{\mathbb E} \def\bbF{\mathbb F} \def\bbG{\mathbb G}
 \def\bbI{\mathbb I} \def\bbJ{\mathbb J} \def\bbK{\mathbb K} \def\bbL{\mathbb L} \def\bbM{\mathbb M} \def\bbN{\mathbb N} \def\bbO{\mathbb O}
 \def\bbP{\mathbb P} \def\bbQ{\mathbb Q} \def\bbS{\mathbb S} \def\bbT{\mathbb T} \def\bbU{\mathbb U} \def\bbV{\mathbb V}
 \def\bbW{\mathbb W} \def\bbX{\mathbb X} \def\bbY{\mathbb Y} \def\kappa{\varkappa} \def\epsilon{\varepsilon}
 \def\phi{\varphi} \def\le{\leqslant} \def\ge{\geqslant}

\def\UU{\bbU}
\def\Mat{\mathrm{Mat}}
\def\tto{\rightrightarrows}

\def\Gr{\mathrm{Gr}}

\def\graph{\mathrm{graph}}

\def\O{\mathrm{O}}

\def\la{\langle}
\def\ra{\rangle}

\def\B{\mathrm B}
\def\Int{\mathrm{Int}}

\def\kre{\epsfbox{hurwitz.6}}
\def\kva{\epsfbox{hurwitz.7}}
\def\tre{\epsfbox{hurwitz.8}}

\def\plusik{\epsfbox{chip.9}}
\def\minusik{\epsfbox{chip.10}}

\def\I{\mathbb I}
\def\M{\mathbb M}
\def\T{\mathbb T}
\def\S{\mathbb S}

\begin{center}
\Large\bf
Infinite symmetric group
and combinatorial descriptions of
 semigroups of double cosets 

\bigskip

\large\sc
Yury A. Neretin\footnote{Supported by the grant FWF, P22122.}
\end{center}

{\small
We show that double cosets of the infinite symmetric group
with respect to some special subgroups admit  natural 
structures of semigroups. We interpret elements of such semigroups
in combinatorial terms (chips, colored graphs,
two-dimensional surfaces with polygonal tiling)
and  describe explicitly the multiplications.
}

\bigskip

The purpose of the paper is to formulate several
facts of the representation theory
of infinite symmetric groups (some references are
\cite{Tho}, \cite{VK}, \cite{Lie}, \cite{Olsh-lieb},
\cite{Olsh-symm}, \cite{Oko}, \cite{KOV}, \cite{Ner-symm})
  in a
natural generality. We consider the infinite symmetric group
$\S_\infty$,
its finite products with itself 
$\S_\infty\times\dots\times \S_\infty$,
and some wreath products. Let $G$, $K$ be groups of this kind.

Sometimes double cosets 
$K\setminus G/K$ admit natural associative product (this is 
a phenomenon, which is usual for infinite dimensional groups%
\footnote{For compact groups double cosets form
 so-called ``hypergroups``, convolution of a pair of double cosets
is a measure on $K\setminus G/K$.}).
We formulate wide sufficient conditions for existence 
of such multiplications.
Our main purpose is to describe such semigroups
 (in fact, categories), in a strange way this can be obtained
in a very wide generality. 

Note that a link between symmetric groups and two-dimensional
surfaces arises at least to Hurwitz \cite{Hur}, see a recent
discussion in \cite{LZ}, see also below Subsection
\ref{ss:belyi}. Categories of polygonal bordisms
were discussed in \cite{Bae}, \cite{Nat}, \cite{Ner-symm}.
Note that they can be regarded of combinatorial analogs 
of  'conformal field theories', see \cite{Seg}, \cite{Ner-book}.

It seems that our constructions can be interesting 
for finite symmetric groups, we give descriptions of various 
double cosets spaces and also produce numerous 'parameterizations'
of symmetric groups.

\section{Introduction}

\COUNTERS

{\bf\punct Double cosets.} Let $G$ be a group,
$K$, $L$ be its subgroup. A {\it double coset} on $G$ is a 
set of the form $KgL$, where $g$ is a fixed element of $G$.
By $K\setminus G/L$ we denote the space of all double cosets.

\sm


{\bf \punct Convolution of double cosets.}
Let $G$ be a Lie group, $K$ a compact subgroup.
Denote by $\cM(K\setminus G/K)$ the space
of all compactly supported
charges (signed measures)
on $G$ invariant with respect to left and right shifts by
to elements of $K$. The space $\cM(K\setminus G/K)$ is an algebra
with respect to the convolution $\mu*\nu$ on $G$. Denote by 
$\Pi\in\cM(K\setminus G/K)$ the probability Haar measure on $K$.
For $\mu\in \cM(K\setminus G/K)$, we have $\Pi*\mu*\Pi=\mu$.

Let $\rho$ be a unitary representation of $G$ in a Hilbert space
$H$. Denote by $H^{[K]}$ the set of all $K$-fixed vectors
in $H$. For a measure $\mu\in \cM(K\setminus G/K)$
consider the operator $\rho(\mu)$ in $H$ given by
$$
\rho(\mu)=\int_G\rho(g)\,d\mu
.$$
By definition,
$$
\rho(\mu)\rho(\nu)=\rho(\mu*\nu).
$$
Next, $\rho(\Pi)$ is the operator $P^{[K]}$ of projection
to the subspace $H^{[K]}$.
Evidently, for $\mu\in \cM(K\setminus G/K)$, we have
$$
\rho(\mu)=\rho(\Pi*\mu*\Pi)=P^{[K]}\rho(\mu)P^{[K]}
.
$$
Therefore, the matrix of $\rho(\mu)$ with respect to the orthogonal decomposition $H=H^{[K]}\oplus( H^{[K]})^\bot$ has the form
$\begin{pmatrix} *&0\\0&0\end{pmatrix}$. In other words, we can regard
the operator $\rho(\mu)$ as an operator
$$
\rho(\mu):H^{[K]}\to H^{[K]}
.
$$
In some special cases ($G$ is a semisimple Lie group, $K$ is 
the maximal compact subgroup; $G$ is a $p$-adic semisimple 
group, $K$ is the Iwahori subgroup, etc.) the algebra 
$\cM(K\setminus G/K)$ and its representation in $H^{[K]}$
are widely explored in representation theory.

\sm


{\bf\punct Noncompact subgroups.}
Now let a subgroup $K$ be noncompact. Then the construction breaks
in both points. First, there is no nontrivial finite $K$-biinvariant
measures on $G$. Second, for faithful unitary representations
of $G$ the subspace $H^{[K]}$ is trivial. 

However, for infinite dimensional groups these objects appear
again  
in another form.


\sm

{\bf\punct Multiplication of double cosets.}
 For some infinite-dimensional groups $G$ 
and for some subgroups $K$ the following facts 
(the {\it ``multiplicativity theorem''})
hold:

\sm

--- There is a natural associative multiplication 
$\frg\circ\frh$ on the set
$K\setminus G/K$.

\sm

--- Let $\frg\in K\setminus G/K$ and $g\in\frg$.
 We define the operator
$$
\rho(\frg):=P^{[K]}\rho(g):\,\, H^{[K]}\to H^{[K]}
$$
Then 
$$
\rho(\frg)\rho(\frh)=\rho(\frg\circ\frh)
.$$

\sm

Multiplicativity theorems were
 firstly observed by Ismagilov
(see \cite{Ism1}, \cite{Ism2}, see also \cite{Olsh-tree})
 for $G=SL(2,Q)$,
$K=\SL(2,Z)$, where $Q$ is complete non-locally
 compact non-archimedian
field, and $Z\subset Q$ is the ring of integer elements.

 Olshanski used semigroups of double cosets 
in the  representation theory of infinite-dimensional classical
groups \cite{Olsh-Opq}--\cite{Olsh-GB} and the representation theory
of infinite symmetric groups in \cite{Olsh-symm},
see also \cite{Olsh-topics}. 
In \cite{Ner-sto}
this approach was extended
to groups of automorphisms of measure spaces.

In \cite{Ner-book} it was given a simple unified proof
of multiplicativity theorems, which covered 
all known cases.
 This also implied that the phenomenon holds under 
quite weak conditions for a pair $G\supset K$.
However it was not clear how to describe double cosets
 $K\setminus G/K$ 
and their products
explicitly. 

Another standpoint of this paper was a  work of Nessonov
\cite{Ness1}--\cite{Ness2}.
For $G=\GL(\infty,\C)\times\dots\times
\GL(\infty,\C)$
(number of factors is arbitrary) with the diagonal subgroup
$K=\U(\infty)$, he classified all $K$-spherical functions on
the group $G$ with respect to $K$.

In \cite{Ner-char} and \cite{Ner-faa} it was described 
multiplication of double cosets for a wide class of pairs
$G\supset K$ of infinite-dimensional classical groups.
In \cite{Ner-symm} parallel
 description was obtained for the pair
$G=\S_\infty\times \S_\infty\times \S_\infty$, $K=\S_{\infty-n}$.

The purpose of this paper is to obtain a similar description
for multiplication of double cosets
for  general pairs related to infinite symmetric groups.


\section{$(G,K)$-pairs and their trains.
 Definitions and a priory theorems}

\COUNTERS

\begin{figure}
 
$$\epsfbox{wreath.1}$$

a) The group $\S_\infty\ltimes (S_k)^\infty$. The group $\S_\infty$
acts by permutations of columns. The normal divisor $(S_k)^\infty$
permute elements in each column. The semi-direct product consists
of permutations preserving partition of the strip into columns

$$\epsfbox{wreath.2}$$

b) The group
$\S_\infty\ltimes (S_{k_1}\times\dots\times S_{k_p})^\infty$.
The
$\S_\infty$ acts by permutations of columns. The normal divisor acts
by permutations inside each sub-column.

$$\epsfbox{wreath.3}$$

c) The set 
$\cup_{i=1}^p\cup_{j=1}^q
\left(\N_i\times \I(\zeta_{ji})\right)$.
Here $q=4$, $p=3$, 
$
Z=\begin{pmatrix}
0&0&1\\
1&5&1\\
2&0&4\\
4&2&1   
  \end{pmatrix}
$

\caption{Reference to Section 2.\label{fig:polosa}}
\end{figure}

\begin{figure}
 $$\epsfbox{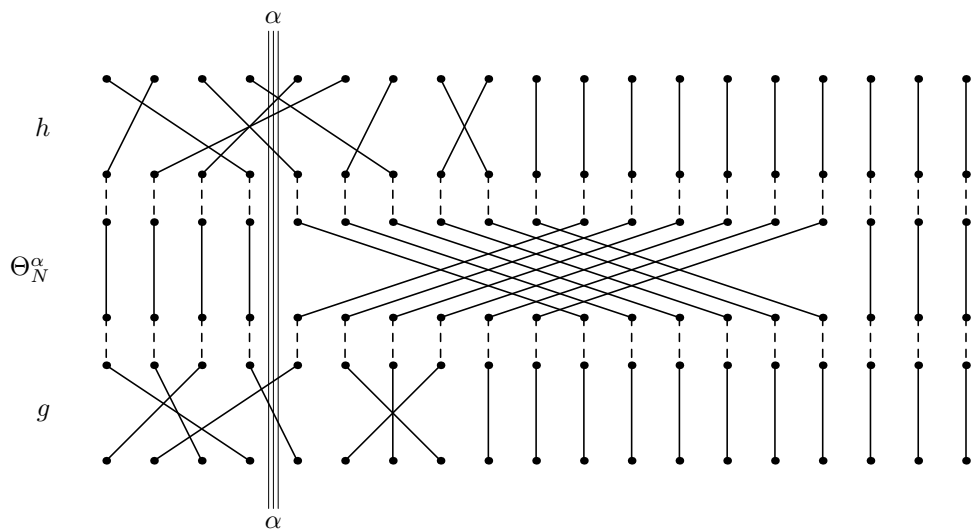}$$

\caption{Reference to Subsection \ref{ss:train}. 
{\it Forcing apart}.
 We draw an element of a symmetric
group as a collection of inclined lines. 
The product is the connection of the corresponding 
 ends. Here $\alpha=\beta=\gamma$.
If $N$ is large, the substitution $\Theta^\alpha_N$ 
sent the set $\supp(h)\cap \{k>\alpha\}$ outside $\supp(g)$
\label{fig:razdvig}.}

\end{figure}

First,  we define some pairs $G\supset K$ of groups,
 for which
multiplicativity theorems
hold.

Secondly, we formulate a priory theorems on the semigroups 
and categories of
double cosets. Proofs are omitted because they are one-to-one
copies of proofs given  in \cite{Ner-faa} for infinite-dimensional
classical groups.

\sm


{\bf \punct Notation.}

\sm

---  $\N$ are natural numbers.

\sm

--- $\M(\alpha)$ is the initial segment $\{1,2,\dots,\alpha\}$
 of $\N$.

\sm

---  $\N_1$, $\N_2$, \dots are disjoint copies
of the set $\N$.

\sm

--- $\M_j(\alpha)$ are initial segments
$\{1,2,\dots,\alpha\}$
 of $\N_j$.
 
\sm

--- $\I(\zeta)$ are  finite sets with $\zeta$ elements.

\sm

--- $\sqcup$, $\coprod$ are symbols for disjoint union of sets. 

\sm

--- $S_n$ is the finite symmetric group.

\sm


{\bf \punct Infinite symmetric group.
\label{ss:s-infty}} 
We denote by $\S_\infty$ the 
group of finite permutations of $\N$. We say that a permutation
$\sigma$ is 
{\it finite}, if $\sigma(i)=i$ for all but finite number of $i$. 
Denote by $\supp(\sigma)$ the {\it support} of $\sigma$,
i.e., the set of $i\in \N$ such that $\sigma(i)\ne i$.

For a countable set $\Omega$ 
denote by $\S_\infty(\Omega)$
 the group of all finite permutations of $\Omega$. 

We represent permutations by infinite $0$-$1$ matrices.

By $\S^\alpha_\infty\subset \S_\infty$
 we denote the group of permutations
having the form
$$
\sigma=\begin{pmatrix}1_\alpha&0\\0&* \end{pmatrix}
,$$
where $1_\alpha$ is the unit $\alpha\times\alpha$ matrix.

Let $K=\S_\infty(\N_1)\times\dots\times \S_\infty(\N_p)$.
For a multi-index $\alpha=(\alpha_1,\dots,\alpha_p)$ we define
the subgroup $K^\alpha\subset K$,
$$
K^\alpha:=
\S_\infty^{\alpha_1}(\N_1)\times\dots\times \S_\infty^{\alpha_p}(\N_p)
.$$


{\bf\punct Topological infinite symmetric group.%
\label{ss:topologic}} We denote by
$\bfS_\infty$ the group of all permutations of $\N$. Define the subgroups 
$\bfS_\infty^\alpha \subset \bfS_\infty$ as above. Define the
topology on $\bfS_\infty$ assuming that $\bfS_\infty^\alpha$ form
a fundamental systems of open neighborhoods of unit. In other words,
a sequence $\sigma_j$ converges to $\sigma$ if for any $k\in \N$
we have $\sigma_j k=\sigma k$ for sufficiently large $j$. The group
$\bfS_\infty$ is a totally disconnected topological group. 

The classification of irreducible unitary representations
of $\bfS_\infty$  was obtained by Lieberman \cite{Lie},
see  expositions in \cite{Olsh-lieb}, 
\cite{Ner-book}.

\begin{theorem}
a) Each irreducible
representation of $\bfS_\infty$ is induced from 
a finite-dimensional representation
$\tau\otimes id$ of a subgroup $S_\alpha\times \bfS^\alpha_\infty$.

\sm

b) Any unitary representation of $\bfS_\infty$ is a direct sum
 of irreducible representations.
\end{theorem}

Note that the quotient space 
$\S_\infty/\S^\alpha_\infty= \bfS_\infty/\bfS^\alpha_\infty$
 is countable, 
and therefore the definition of induced representations
 survives (see, e.g.,
\cite{Kir}, 13.2).

In a certain sense, the Lieberman theorem opens and closes
the representation theory of the group $\bfS_\infty$. However it 
is an important element of wider theories.

\sm


{\bf\punct Reformulations of continuity.%
\label{ss:reformulation-continuity}}
Let $\rho$ be a unitary representation of $\S_\infty$
in a Hilbert space $H$. Denote by $H^\alpha\subset H$
the subspace of all 
$\S^\alpha_\infty$-fixed vectors. We say that a representation
$\rho$ is 
{\it admissible} if $\cup H^\alpha$ is dense in $H$.

Denote by $\bfB_\infty$ the semigroup of matrices
 composed of $0$ 
and $1$
such that each row and each column contains $\le 1$ units. 
We equip $\bfB_\infty$ with the topology of element-wise convergence,
the group $\S_\infty$ is dense in $\bfB_\infty$.

\begin{theorem}
\label{th:admissibility}
{\rm(see \cite{Olsh-lieb}, \cite{Ner-book})}
 The following conditions are equivalent:

\sm

1) $\rho$ is continuous in the topology of
$\bfS_\infty$;

\sm

2) $\rho$ is admissible;

\sm

3) $\rho$ admits a continuous extension to the semigroup
$\bfB_\infty$.
\end{theorem}

This statement has a straightforward extension to products
of symmetric groups
$\bfK=\bfS_\infty\times\dots\times \bfS_\infty$

\sm


{\bf\punct  Wreath products.%
\label{ss:wreath}}
Let $U$ be a finite group.
Consider the countable  {\it direct
product} 
$\mathbf{U^{\mathbf\infty}}:=U\times U\times U\dots,$
it is the group, whose elements are infinite sequences
$(u_1,u_2,\dots)$.
Consider also the {\it restricted product} $U^\infty$, 
whose elements are
sequences such that $u_j=1$ starting some place.
The group $\mathbf{U^{\infty}}$ is equipped with the topology
of direct product, the group ${U^{\infty}}$ is discrete.

Permutations of sequences $(u_1,u_2,\dots)$
induce automorphisms of  $U^\infty$ and $\mathbf{U^\infty}$.
Consider the semidirect products
$K:=\S_\infty\ltimes U^\infty$ and $\bfK:=\bfS_\infty
\ltimes\mathbf{U^\infty}$, see, e.g., \cite{Kir}, 2.4,
they are called {\it wreath products} of $\S_\infty$ and $U$.

\sm

{\sc Example.} The  infinite {\it hyperoctahedral group}
 is a wreath product of $\S_\infty$ and $\Z_2$.
\hfill $\square$

\sm

Our main example is the wreath product of $\S_\infty$ and
a finite symmetric group  $S_k$.  We realize it
as a group of finite permutations of the
$\N\times\I(k)$,
 The group $\S_\infty$ acts
by permutations of $\N$, and subgroups $S_k\subset (S_k)^\infty$
 act
by permutations of sets $\{m\}\times\{1,\dots,k\}$,
see Figure \ref{fig:polosa}.a.

\sm


{\bf\punct Representations of wreath products.}
For $K=\S_\infty\ltimes U^\infty$ we define a subgroup
$K^\alpha$ 
 as the semidirect product of
$\S^\alpha_\infty$ and the subgroup 
$$
U^{\infty-\alpha}:=
\underbrace{1\times\dots 1}_{\text{$\alpha$ times}}
\times U\times U\times \dots
.
$$

We define {\it admissible representations} as above.
Again, a unitary representation of $K=\S_\infty\ltimes U^\infty$ 
is admissible if and only if it 
is continuous in the topology of the group
$\bfK=\mathbf{S_\infty\ltimes U^\infty}$.


\sm

{\bf\punct  $(G,K)$-pairs,%
\label{ss:G-K}} see Figure \ref{fig:polosa}.
Fix positive integers $q$, $p$
 and $q\times p$-matrix
$$
Z:=\{\zeta_{ji}\}
$$
consisting of non-negative integers. Assume that
it has no zero columns and no zero rows.

Fix a collection 
$\Lambda=
\begin{pmatrix}\lambda_1\\ \vdots\\ \lambda_q\end{pmatrix}$
of nonnegative integers.
Fix  sets $L_1$,\dots, $L_q$,
such that $L_j$ has $\lambda_j$ elements.
Consider the collection of  sets
$$
\N_i\times \I(\zeta_{ij})
.$$
Assume that $\Omega_j$ is the disjoint union 
$$
\Omega_j:= L_j\sqcup \coprod_{i\le p}
 \left( \N_i\times \I(\zeta_{ji})
\right)
$$
and assume that all sets $\Omega_j$ obtained
in this way are mutually disjoint.
Set
$$
G:=G[Z,\Lambda]=\prod_{j=1}^q \S_\infty(\Omega_j)
.
$$
Next, set
$$
K^\circ=K^\circ[Z]:=\prod_{i=1}^p\left(
\S_\infty(\N_i) 
\ltimes \Bigl(\prod_j S_{\zeta_{ij}}\Bigr)^\infty\right)
.
$$
We have a tautological embedding $K^\circ\to G$.

Also, we set
$$
K^\circledast[Z]:= \prod_{i=1}^p \S_\infty(\N_i)\subset K^\circ[Z]
.$$

Below $(G,K)$ denotes a pair
(group, subgroup) of the form
$$
(G,K)=(G[Z,\lambda],K^\circ[Z]) \qquad
\text{or}
\qquad 
(G,K)=
(G[Z,\lambda],K^\circledast[Z])
.$$

{\sc Remark.}
It is also possible to consider intermediate
wreath products $K$ between  $K^\circ[Z]$ and
$K^\circledast[Z]$, below we consider one example from
 this
zoo.

\sm


{\bf \punct Colors, smells, melodies.%
\label{ss:melody}} We wish to draw figures, also we 
want to have  more flexible language.

\sm

a) We assign to each $\Omega_j$ a {\it color},
say, red, blue, white, red, green, etc. We also think
that a color is an attribute of all points of $\Omega_j$.
We denote colors by $\gimel_j$.

\sm

b) Next, we assign to each $\N_i$ a {\it smell}
$\aleph_i$,
say, Magnolia, Matricana, Pinus, Ledum, Rafflesia, etc.  
On figures we denote smells by $\tre$, $\kre$, $\kva$, \dots.
We also thick that a smell $\aleph_i$ 
is an attribute of all points of 
$\left(\N_i\times \I(\zeta_{ji})\right)\subset \Omega_j$.

\sm

c) Orbits of a group $\S_\infty(\N_i)$ on $\Omega_j$
are one-point orbits or countable homogeneous
spaces $\S_\infty/\S^1_\infty\simeq\N$.
We assign to each countable orbit a {\it melody},
say,  violin, harp, tomtom, flute, drum,\dots.
On figures we draw melodies by symbols $\heartsuit$,
$\succ$, $\nabla$, $\sharp$, $\ddag$, etc.
Note that a melody makes sense after fixation 
of a smell and a color.

\sm

{\sc Example.} See Figure \ref{fig:polosa}.c.
We have a $4\times 3$ table. 
Rows are distinguished by colors, columns are 
separated by smells. Rows inside each box
are numerated by melodies. \hfill$\square$


\sm

{\bf\punct Admissible representations.%
\label{ss:admissible-GK}}
Let $\rho$ be a unitary representation of $G$. We say that
$\rho$ is a {\it $K$-admissible representation} 
 if the restriction of $\rho$ to $K$ is admissible.
Equivalently, we say that {\it $\rho$ is a representation of the pair
$(G,K)$.}

\sm


{\bf\punct Reformulation of admissibility in terms
of continuity.%
\label{ss:continuity-GK}}
The embedding $K\to G$ admits an extension to the map
 $\bfK\to\bfG$
of the corresponding completions.
Consider the group 
$G\cdot \bfK$ generated by
$G$ and $\bfK$,
$$
G\subset G\cdot \bfK\subset \bfG
.$$
 Any element of $G\cdot \bfK$ admits a (non-unique) 
representation as
$g\bfk$, where $g\in G$, $\bfk\in K$.

We consider the natural topology 
on the subgroup $\bfK$ and assume that $\bfK$ is an
 open-closed subgroup
in $G\cdot \bfK$.

\begin{proposition}
A unitary representation of $G$ is $K$-admissible if and only 
if it is continuous in the above sense.
\end{proposition}

{\sc Proof.} Let $\rho$ be an admissible representation of $G$ 
in a Hilbert space $H$. We define the action 
of $\bfK$ as a continuous extension of the action of $K$.
\hfill $\square$


\sm

{\bf\punct Lemma on admissibility.%
\label{ss:l-admissibility}} 

\begin{lemma}
Let $\rho$ be an irreducible unitary representation of $G$.
If some $H^\alpha\ne 0$, then the representation $\rho$
is $K$-admissible.
\end{lemma}

{\sc Proof.} Consider the subspace $\cH:=\cup H^\alpha$. 
Fix $g\in G$.
For sufficiently large $\beta$, the element $g$ commutes
with $K^\beta$. Therefore $\cH$ is $g$-invariant.
The closure of $\cH$ is a subrepresentation.%
\hfill $\square$

\begin{corollary}
If an irreducible unitary representation of $G$ has a
$K$-fixed vector, then it is $K$-admissible.
\end{corollary}


{\bf\punct Existing representation theory.%
\label{ss:existing}}
 The most interesting object
of existing theory is the pair 
 $G=\S_{\infty}\times \S_{\infty}$, $K=\S_\infty$ 
is the diagonal subgroup, \cite{Olsh-symm}, \cite{Oko}, \cite{KOV}).
In our notation 
$Z=\begin{pmatrix}   1\\1   \end{pmatrix}$.
The representation theory of this pair includes 
also earlier works
on Thoma characters (see \cite{Tho}, \cite{VK}, \cite{VK2}).

Olshanski \cite{Olsh-symm} also considered pairs:

\sm

--- $G=\S_{\infty+1}\times \S_{\infty}$, $K=\S_\infty$;
in our notation, $Z=\begin{pmatrix}   1\\1   \end{pmatrix}$,
 $\Lambda=\begin{pmatrix}   1\\0 \end{pmatrix}$.

\sm

--- $G=\S_{2\infty}$, $K=\S_\infty\times \S_\infty$; in our
notation $Z=\begin{pmatrix}   1&1   \end{pmatrix}$.

\sm

--- $G=\S_{2\infty}$, 
$K=\S_\infty\ltimes\Z_2^\infty$ and also $G=\S_{2\infty+1}$
with the same subgroup $K$. In our notation, $Z=(2)$
and $\Lambda=0$ or $1$.

\sm

{\it In all these cases, pairs $(G,K)$ are  limits 
of spherical pairs
of finite groups.}

\sm

The author in \cite{Ner-symm} considered 
the case $G=\S_\infty\times\dots\times \S_\infty$ with the diagonal
subgroup $K=\S_\infty$. In our notation,
$Z=\begin{pmatrix}   1\\\vdots\\ 1   \end{pmatrix}$.


\sm

{\bf\punct Train.%
\label{ss:train}}
Consider (see Figure \ref{fig:razdvig})
 the following matrix $\Theta^{[\alpha]}_N\in \S_\infty$
of the size
$(\alpha+N+N+\infty)\times (\alpha+N+N+\infty)$,
$$
\Theta^{[\alpha]}_N=
\begin{pmatrix}
1_\alpha&0&0&0
\\
0&0&1_N&0
\\
0&1_N&0&0
\\
0&0&0&1_\infty
\end{pmatrix}\in \S_\infty
.
$$
In fact, $\Theta^{[\alpha]}_N$ is contained in $K^\alpha$.

Consider a pair $(G,K)$.
For a multi-index $\alpha=(\alpha_1,\dots,\alpha_p)$ we denote by
$\Theta^{[\alpha]}_N$ the element
$$
\Theta^{[\alpha]}_N=
\left(\Theta^{[\alpha_1]}_N,\dots,\Theta^{[\alpha_p]}_N\right)\in K
$$
Again, $\Theta^{[\alpha]}_N\in K^\alpha$.

Fix multi-indices  $\alpha$, $\beta$, $\gamma$. Consider double cosets
$$
\frh\in K^\beta\setminus G/K^\alpha,
\qquad
\frg\in K^\gamma\setminus G/K^\beta
,$$
and choose their representatives
$g\in\frg$,  $h\in\frh$.
Consider the sequence
$$
f_N=g\Theta^{[\alpha]}_Nh\in G
.$$
Consider the double coset $\frf_N$ containing $f_N$,
$$
\frf_N\in K^\gamma\setminus G/K^\alpha
.$$

\begin{theorem}
\label{th:product}
{\rm a)} The sequence $\frf_N$ is eventually constant.

\sm

{\rm b)} The limit 
$$\frg\circ\frh:=\lim_{N\to\infty}\frf_N$$
does not depend on a choice of representatives $g$, $h$.

\sm

{\rm c)} The product $\circ$ obtained in this way is associative,
i.e., for any
$$
\frg\in K^\delta\setminus G/K^\gamma,
\qquad
\frh\in K^\gamma\setminus G/K^\beta,
\qquad
\frf\in K^\beta\setminus G/K^\alpha,
,$$
we have
$$
(\frg\circ \frh) \circ \frf = \frg\circ (\frh \circ \frf)
.$$
\end{theorem}

Thus we obtain a category $\T(G,K)$, whose objects are multiindices
$\alpha$ and morphisms $\alpha\to\beta$ are elements
of $K^\beta\setminus G/K^\alpha$.
We say that $\T(G,K)$ is the {\it train} of the pair
$(G,K)$.

\sm


{\bf\punct Involution in the train.%
\label{ss:involution}}
The map $g\mapsto g^{-1}$ induces a map of quotient spaces
$
K^\alpha\setminus G/K^\beta
\to K^\beta\setminus G/K^\alpha
$, we denote it by
$$
\frg\mapsto \frg^\square
.$$
Evidently,
$$
(\frg\circ \frh)^\square=\frh^\square \frg^\square
.$$


\sm

{\bf\punct Representations  of the train.%
\label{ss:repres-train}}
Now let $\rho$ be a unitary representation of the pair $(G,K)$.
We define subspaces $H^\alpha$ as above, denote by $P^\alpha$
the operator of orthogonal projection to $H^\alpha$.
For
$\frg\in K^\beta\setminus G/K^\alpha$,  choose its representative
$g\in\frg$. 
Consider
the operator
$$
\ov\rho_{\alpha,\beta}(g):=P^\beta\rho(g):
\,\,
H^\alpha\to H^\beta
.
$$
By definition, we have
\begin{equation}
\|\ov\rho_{\alpha,\beta}(g)\|\le 1
\label{eq:contractive}
.
\end{equation}

\begin{theorem}
\label{th:representation}
{\rm a)} An operator $\ov\rho_{\alpha,\beta}(g)$
 depends only on a double coset
$\frg$ containing $g$.

\sm

{\rm b)} We get a representation of the  category $\T(G,K)$,
i.,e., for any
$$
\frg\in K^\gamma\setminus G/K^\beta,
\qquad
\frh\in K^\beta\setminus G/K^\alpha
$$
the following identity holds
$$
\ov\rho_{\beta,\gamma}(\frg)
\,\ov\rho_{\alpha,\beta}(\frh)=
\ov \rho_{\alpha,\gamma}(\frg\circ\frh)
.
$$

{\rm c)} We get a $*$-representation. i.e.,
$$
\bigl(\ov\rho_{\alpha,\beta}(\frg)\bigr)^*=
\ov\rho_{\beta,\alpha}(\frg^\square)
.
$$

{\rm d)} $\rho(\Theta_N^\alpha)$
weakly converges to the projection $P^\alpha$.
\end{theorem}

\begin{theorem}
\label{th:bijection}
Our construction provides a bijection between
the  set of all unitary representation of the pair $(G,K)$ 
and the set of all $*$-representations
of the category $\T(G,K)$ satisfying the condition
{\rm (\ref{eq:contractive}).}
\end{theorem}

We omit proofs of Theorems \ref{th:product}--\ref{th:bijection}
and Theorem \ref{th:spherical} formulated below, because
 proofs are literal copies of proofs in \cite{Ner-faa}.

Our main purpose is to give explicit description of trains,
also we give some constructions of representations of groups.


\sm

{\bf\punct Sphericity.%
\label{ss:sphericity}} 

\begin{theorem}
\label{th:spherical}
Consider a pair $(G,K)=(G[Z,\Lambda],K^\circ[Z])$,
or $(G[Z,\Lambda],K^\circledast[Z])$ as above. 
If $\Lambda=0$, then the pair $(G,K)$ is spherical.
In other words, for any irreducible unitary representation of
$(G,K)$ the dimension of the space of $K$-fixed
vectors is $\le 1$.
\end{theorem}

\sm



\sm

{\bf\punct Further structure of the paper.}
In Sections 8, 9 we present the description
 of  trains for arbitrary pairs
 $$(G(Z,\Lambda),K^\circ[Z]),\qquad
(G(Z,\Lambda),K^\circledast[Z]).
$$
We also present examples of representations
and spherical functions.

I am afraid that it is difficult to understand constructions
in such generality. For this reason,
in Sections 3--7 we consider simple special cases: 
pairs $(G,K)$ connected
with matrices
$$
Z=\begin{pmatrix}1&\dots&1 \end{pmatrix},
\qquad
Z=\begin{pmatrix}2&\dots&2 \end{pmatrix},
\qquad
Z=\begin{pmatrix}3\end{pmatrix}
.$$

Well-representative examples are contained in Section 7.

\sm

Now we give 4 examples (details are below), to outline
the language used for description of double coset spaces.

\sm

{\sc Example 1.} Take collection of beads
of $n$ colors. 
 We say that a chaplet is a cyclic sequence of $2k$ beads,
we also equip beads of a chaplets by interlacing signs $\pm$.
Consider the set $\Xi$, whose elements are (non-ordered)
 collections of chaplets. Then
$\Xi\simeq K^\circ[Z]\setminus G[Z,0]/K^\circ[Z]$,
where $Z=\begin{pmatrix}2&\dots&2\end{pmatrix}$.

\sm

{\sc Example 2.} Consider the set $\Delta$ 
of (generally, disconnected) oriented triangulated
closed two-dimensional surfaces
with additional data: triangles are colored in black and white
checker-wise. We also forbid a triangulation 
of a sphere into two triangles. 
Let $(Z)=3$, $\Lambda=0$.
We consider the subgroup
$K^\odot:=\S_\infty\ltimes (\Z_3)^\infty$
intermediate between $K^\circ[Z]$
and $K^\circledast[Z]$.
The set $\Delta$
is in a canonical one-to-one correspondence with
the double coset space
$K^\odot\setminus G[Z,0]/K^\odot$.

\sm

{\sc Example 3.} Consider the set $\Phi$ of all 
tri-valent graphs, whose vertices are colored black and white
(and neighboring vertices have different colors). We forbid
components having only two vertices. The set
$\Phi$ is in one-to-one
correspondence
with $\Xi\simeq K^\circ[Z]\setminus G[Z,0]/K^\circ[Z]$, where
$Z=(3)$.

\sm

{\sc Example 4.} (Belyi data) Consider the Riemannian sphere
$\ov\C=\C\cup\infty$
with 3 distinguished points $0$, $1$, $\infty$.
Consider the set
$\Psi$ all (generally, disconnected)
 ramified coverings $\Gamma\to\ov \C$
with 3 branching points $0$, $1$, $\infty$
(requiring that the covering is nontrivial on each component
of $\Gamma$).
Then $\Psi$ is in a canonical one-to-one correspondence
with the set $K^\circledast[Z]\setminus G[Z,0]/K^\circledast[Z]$,
where $Z:=\begin{pmatrix}
        1\\1\\1       \end{pmatrix}$.


\sm

{\bf \punct Remarks on generality.}
Note, that our generality is not maximal.

\sm

1. Our constructions can be easily extended to the case
$p=\infty$ or $q=\infty$.

\sm

2. However, it seems
that number of melodies of a fixed color
and a fixed smell must be finite. Also, in
our picture,
 orbits of $\S_\infty(\N_i)$ on $\Omega_j$
are fixed points or $\S_\infty(\N_i)/\S_\infty^1(\N_i)$.
It seems that this is necessary, otherwise
we have embeddings
$\S_\infty(\N_i)\to\bfS_\infty(\Omega_j)$ and not to 
$\S_\infty(\Omega_j)$.

\sm

3. Let $L$, $M\subset G[Z,\Lambda]$ be subgroups of the type
described above and $L$, $M\supset K^\circledast[Z]$. Then
double cosets $L\setminus G[Z,\Lambda]/M$ admit a description
on the language of Section 8.


\section{First example. Symmetric group and Young subgroup}

\COUNTERS

\begin{figure}
$$
\epsfbox{hurwitz.1}
$$

a) We present two copies of the set $\Omega$.
An element $g\in \S_\infty(\Omega)$ is drawn as a collection
of segments connecting $\omega'$ and $(g\omega)^{''}$.
The set $\Omega'$ (respectively $\Omega^{''}$) is split 
into the pieces in a natural way. The first piece
is $\Omega'_{[\alpha]}$, any  point of the remain
 is contained in some set $\N_i$ 
and has the corresponding smell.
We draw smells as $
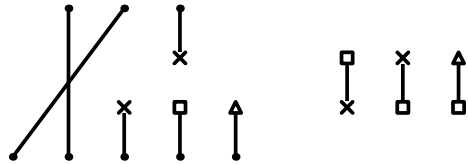
\epsfbox{hurwitz.6}$, 
$\epsfbox{hurwitz.7}$, $\epsfbox{hurwitz.8}$.

$$
\epsfbox{hurwitz.2}
$$

b) We forget part of information about segments.
We remember 

1. ends that are contained in 
$\Omega'_{[\alpha]}$, $\Omega''_{[\beta]}$.

2. smells of remaining ends.

3. direction up/down of segments.
 
Segments whose ends have the same
smell are removed.

$$
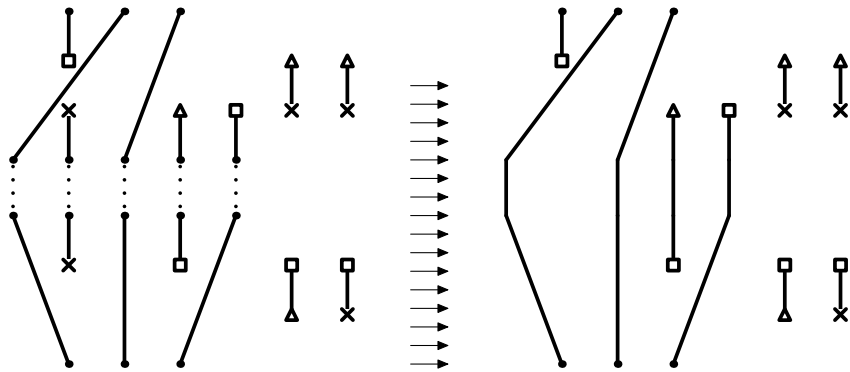
\epsfbox{hurwitz.3}\qquad \epsfbox{hurwitz.5}\qquad \epsfbox{hurwitz.4} 
$$

c) Two diagrams (left figure) and their product (right figure).

\caption{Reference to Section 3.\label{fig:1}}

\end{figure}

{\bf \punct Group.%
\label{ss:group-ness}} Now $q=1$ and $p$ is arbitrary,
$$
Z=\begin{pmatrix}
    1&\dots&1
   \end{pmatrix},
\qquad \text{$\lambda\ge 0$ is arbitrary}
.$$
Therefore, 
 $\Omega=\N_1\sqcup\dots\sqcup \N_p\sqcup L$,
and
\begin{multline*}
(G,K)=(G[Z,\Lambda], K^\circ[Z])=
(G[Z,\Lambda], K^\circledast[Z])
=\\=
\bigl(\S_\infty(\Omega),\,
 \S_\infty(\N_1)\times\dots\times \S_\infty(\N_p)
\bigr).
\end{multline*}
 We assign a smell $\aleph_i$
to each $\N_i$. On some figures we draw symbols
$\epsfbox{hurwitz.6}$, $\epsfbox{hurwitz.7}$, 
$\epsfbox{hurwitz.8}$
instead of smells.

\sm


{\bf\punct Description of the category.%
\label{ss:category-young}}
For a subgroup $K^\alpha\subset K$
denote by
$\Omega_{[\alpha]}\subset \Omega$ 
 the subset fixed by $K^\alpha$,
$$
\Omega_{[\alpha]}=L\sqcup \M_1(\alpha_1)\sqcup
\dots \sqcup \M_p(\alpha_p)
.
$$

Consider two copies $\Omega'$ and $\Omega{''}$ 
of the set $\Omega$. 
 For  a point of $\omega\in \Omega$
we denote its copies by $\omega'$ and $\omega^{''}$.

For any $g\in \S_\infty(\Omega)$
we construct an oriented one-dimensional manifold
$\Xi_{\alpha,\beta}[g]$
 with 
boundary
(in fact, $\Xi_{\alpha,\beta}[g]$ is a union of segments).
It is easier to look to Figure \ref{fig:1}, but we present
 a formal description.

\sm

a) If  $\omega\in \Omega_{[\alpha]}$, 
$g\omega \in \Omega_{[\beta]}$,
we draw a segment, whose origin is $\omega'$ 
and the end is $(g\omega)^{''}$.

\sm

b) If 
$\omega\in \Omega_{[\alpha]}$, 
$g\omega \notin \Omega_{[\beta]}$, then we draw a segment starting
in $\omega'$ and mark another end of the segment by the smell of
$g\omega$.

\sm

c) Let $\omega\notin \Omega_{[\alpha]}$, 
$g\omega \in \Omega_{[\beta]}$. 
Then we draw a segment, mark its origin by the smell of
$\omega$, another end of the segment is $(g\omega)^{''}$.

\sm

d) Let
$\omega\notin \Omega_{[\alpha]}$, 
$g\omega \notin \Omega_{[\beta]}$.
If smells of $\omega$ and $g\omega$ coincide, 
then we do not draw a segment.

\sm

e)  Let
$\omega\notin \Omega_{[\alpha]}$, 
$g\omega \notin \Omega_{[\beta]}$.
If smells of $\omega$ and $g\omega$ are different,
we draw a segment whose origin has the smell
of $\omega$ and the  end has the smell of $g\omega$. 

\sm

We say that points of $\Omega_{[\alpha]}'$ are {\it entries}
of $\Xi_{\alpha,\beta}(g)$
and points of $\Omega_{[\beta]}^{''}$ are {\it exits}.

Due d) we get only finite collection
of segments.

\begin{lemma}
If $g_1$, $g_2$ are contained in the same element
of
$K^\beta\setminus G/K^\alpha$,
then $\Xi_{\alpha,\beta}[g_1]=\Xi_{\alpha,\beta}[g_2]$.
\end{lemma}

This is obvious. \hfill $\square$

\sm

Now consider two elements
$\frg\in K^\beta\setminus G/K^\alpha$, 
$\frh\in K^\gamma\setminus G/K^\beta$.
To obtain the diagram corresponding
$\frh\circ\frg$, we identify
each exit of $\Xi_{\alpha,\beta}(\frg)$
with the corresponding  entry  of
$\Xi_{\beta,\gamma}(\frg)$ and get a new 
oriented manifold. It remains to remove segments
whose origin and end have the same smell, see 
Fig.\ref{fig:1}.c).

\begin{proposition}
This product is the product in the train
$\T(G,K)$ of the pair $(G,K)$.
\end{proposition}

{\sc Proof.} This is obvious, see Fig. \ref{fig:razdvig}.
\hfill $\square$


\sm

{\bf\punct  The involution in the category.}
We change entries and exits, and reverse  orientations
of the segments.


\section{Tensor products of Hilbert spaces}

\COUNTERS

{\bf\punct Definition of tensor products.%
\label{ss:tensor-products}}
Let $H_1$, $H_2$, \dots be a countable collection
of Hilbert spaces (they can be finite-dimensional or 
infinite-dimensional). Fix a unit vector $\xi_k\in H_k$
in each space.
The tensor product
$$
(H_1,\xi_1)\otimes(H_2,\xi_2)\otimes(H_3,\xi_3)\otimes
\dots
$$
is defined in the following way.
We choose an orthonormal basis $e_j[k]$
in each $H_k$, assuming $e_1[k]=\xi_k$.
Next, we consider the Hilbert space with orthonormal basis
$$
e_{\alpha_1}[1]\otimes e_{\alpha_2}[2]\otimes\dots
$$
such that $e_{\alpha_N}[N]=\xi_N$ for sufficiently
large $N$ (note, that this basis is countable).

Distinguished vectors are necessary for the definition 
(otherwise we get a non-separable object).

Construction essentially depends on a choice of
distinguished vectors.
The spaces
$
\otimes(H_k,\xi_k)
$
and
$
\otimes(H_k,\eta_k)
$
are canonically isomorphic if and only if 
$$
\sum|\la \xi_j,\eta_j\ra-1|<\infty
.$$
In particular, we can omit distinguished vectors
in a finite number of factors (more precisely we
can choose them in arbitrary way).

All subsequent definitions are standard
(see the  paper of von Neumann \cite{vN} with watching
of all details or short introduction in \cite{Gui}).

\sm 


{\bf\punct Action of symmetric groups in tensor products.}
The symmetric groups $S_n$ act in tensor powers $H^{\otimes n}$
 by permutations of
factors. This phenomenon has a straightforward analog.

We denote by 
$$(H,\xi)^{\otimes\infty}:=
(H,\xi)\otimes (H,\xi)\otimes\dots
$$
the infinite symmetric power of $(H,\xi)$.

\begin{proposition}
\label{pr:tensor-1}
{\rm a)}
The complete symmetric group
$\bfS_\infty$ acts in
$ (H,\xi)^{\otimes\infty}$ by permutation
of factors. The representation is continuous 
with respect to the topology of $\bfS_\infty$.

\sm

{\rm b)} The vector $\xi^{\otimes\infty}$ is a unique
$\bfS_\infty$-fixed vector in $ (H,\xi)^{\otimes\infty}$.

\sm

{\rm c)} The subspace of $\bfS^\alpha_\infty$-fixed vectors is 
$$
H^{\otimes \alpha} \otimes \xi^{\otimes\infty}
.$$
\end{proposition}

\begin{proposition}
\label{pr:tensor-2}
Fix a sequence $\xi_k$ of unit vectors  in 
a Hilbert space 
$H$. The  symmetric group $\bbS_\infty$
acts in the tensor product
$\otimes_k (H,\xi_k)$ by permutations 
of factors.
\end{proposition}

Emphasis that in this case there is no action
of the complete symmetric group $\bfS_\infty$.


\sm

{\bf\punct First application to
 $(G,K)$-pairs.%
\label{ss:rep-young}}
Now consider the same objects $\Omega$, $G$, $K$,
 as in the previous section, recall that
$$
\Omega=L\sqcup\N_1\sqcup\dots\sqcup\N_p
.$$

Consider a  Euclidean space
$H$ of dimension $\le p$ and a  collection of unit vectors 
$\eta_1$, \dots $\eta_p\in H$
generating $H$. 
Consider a countable collection of copies
$H_\omega$
of $H$ enumerated by elements of the set $\Omega$.
For each $\omega\in\Omega$ we choose an element
$\xi_\omega$ according the following rule:

\sm

--- if $\omega\in \N_i$, then $\xi_\omega=\eta_i$,

\sm

--- if $\omega\in L$, then we choose
a unit vector $\xi_\omega$ in an arbitrary way.

\sm

Consider the tensor product
$$
\cH:=
\bigotimes_{\omega\in\Omega} (H,\xi_\omega)=
\bigotimes_{\omega\in L} H
\otimes \bigotimes_{i=1}^p
(H,\eta_i)^{\otimes\infty}
.
$$

The group $G=\S_\infty(\Omega)$ acts in 
$\cH$ by permutation of factors
by Proposition \ref{pr:tensor-2}.

Each group $\bbS_\infty(\N_i)$
acts in $\cH$ by permutations
of factors in
$$
(H,\eta_i)^{\otimes \infty}
$$
Therefore, we get the action of $\bfK=\prod_{i=1}^p \bfS_\infty(\N_i)$
in $\cH$.

Thus, the  group $G\cdot\bfK$ 
 acts in the tensor product $\cH$
by permutations of factors. 

The parameter determining such representations is the  Gram matrix
$$
a_{kl}=\la\eta_k,\eta_l\ra
$$
of the system $\eta_i$.

\sm


{\bf\punct  Super-tensor products}
(a subsection for experts in representation theory).
Below we construct representation for arbitrary $(G,K)$-pairs
related to symmetric groups. All constructions
described below can be extended to super-tensor products
as it is explained in \cite{Olsh-symm}, \cite{Ner-symm}.

For the $(G,K)$-pair discussed now this gives nothing.


\section{Chips}

\COUNTERS

Construction of this section is an extension
of Olshanski \cite{Olsh-symm}. We use the term ``chip'' following
Kerov \cite{Ker}, note that the first version of chips%
\footnote{A matching on the union of two sets.} was
introduced by R. Brauer in \cite{Bra}.

\sm


{\bf \punct The group.\label{ss:chips1}}
Take
$Z=\begin{pmatrix}2&\dots&2\end{pmatrix}
$, let $\Lambda$ be arbitrary.
 Consider the pair 
$(G,K)=(G[Z,\Lambda]), K^\circ[Z])$.
Now 
$$
\Omega =L\sqcup \coprod_{i=1}^p (\N_i\times \I(2))
$$
and $K$ is a product of hyperoctohedral groups.
$$
K=\prod_{i=1}^p
\left(\S_\infty(\N_i)\ltimes \Z_2^\infty\right)
.$$

By $\omega\mapsto\omega^\circ$ we denote the natural involution
in $\Omega$, it acts by the transposition of two elements
of $\I(2)$, this involution fixes points of the set $L$. 

 For each 
$\N_i$,  we attribute a smell $\aleph_i$,
say $\epsfbox{hurwitz.6}$,  $\epsfbox{hurwitz.7}$,
$\epsfbox{hurwitz.8}$, etc.


\sm

{\bf\punct Description of the train.%
\label{ss:chips-train}}
For a subgroup $K^\alpha$ denote by $\Omega_{[\alpha]}$ the set
of all point fixed by $K^\alpha$,
$$
\Omega_{[\alpha]}:=
L\sqcup\coprod_{i=1}^p 
\bigl(\M_i(\alpha_i)\times \I(2)\bigr)
.
$$
Such sets are objects of the train $\T(G,K)$.

We wish to assign a graph with additional data
to any double coset in $K^\beta\setminus G/K_\alpha$.
The construction is explained on Figure \ref{fig:3}.
However, we repeat this formally.

Consider two copies of $\Omega$, say $\Omega'$ and $\Omega^{''}$.
Consider the 'set of entries' $\Omega_{[\alpha]}'$ and
the 'set of exits' $\Omega_{[\beta]}^{''}$.


\begin{figure}
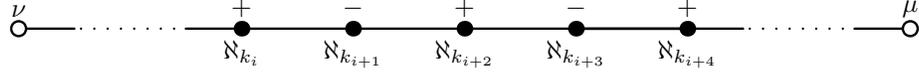
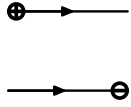
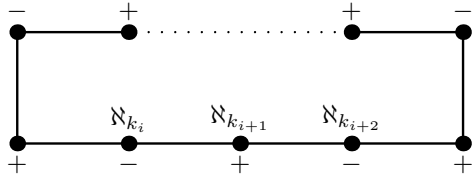
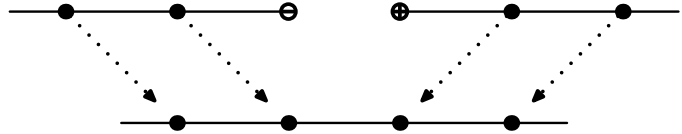

$$
\epsfbox{chip.4}
$$

a) An open chain.

$$
\epsfbox{chip.5}
$$

$$
\epsfbox{chip.6}
$$

b) Ends of an open chain.

$$
\epsfbox{chip.7}
$$

c) A closed chain.

$$
\epsfbox{chip.8}
$$

d) Gluing of two chains.

\caption{Reference to Section 5.\label{fig:2}}
\end{figure}

\sm

a) A {\it closed chain} is a finite cyclic sequence
 of distinguished  points (vertices) and
 segments connecting neighboring vertices.
A vertex has two labels: a smell $\aleph_{k_j}$
and $\pm$.
Symbols $+$ and $-$  interlace
(therefore  a length of a closed chain is even
 and there are only
two possible choices of signs).

\sm

b) An {\it open chain} has similar properties. Interior 
vertices are labeled by smells
and $\pm$.
 Ends are points of the set
$\Omega_{[\alpha]}'\cup \Omega_{[\beta]}^{''}$.
If an end is an 'entry', i.e., an element of 
$\Omega_{[\alpha]}'$, we mark it by $\epsfbox{chip.9}$.
If an end is an exit, then we mark it by $\epsfbox{chip.10}$.
Pluses and minuses interlace. 

If both ends are `entries' (or both are 'exits'),
 then the length 
of the chain is even. If one of the ends is an
 entry and another one is an exit, then 
the length is odd.
 Note that there is a unique possible choice of signs if
we fixed types of ends.

\sm

\begin{figure}
$$
\epsfbox{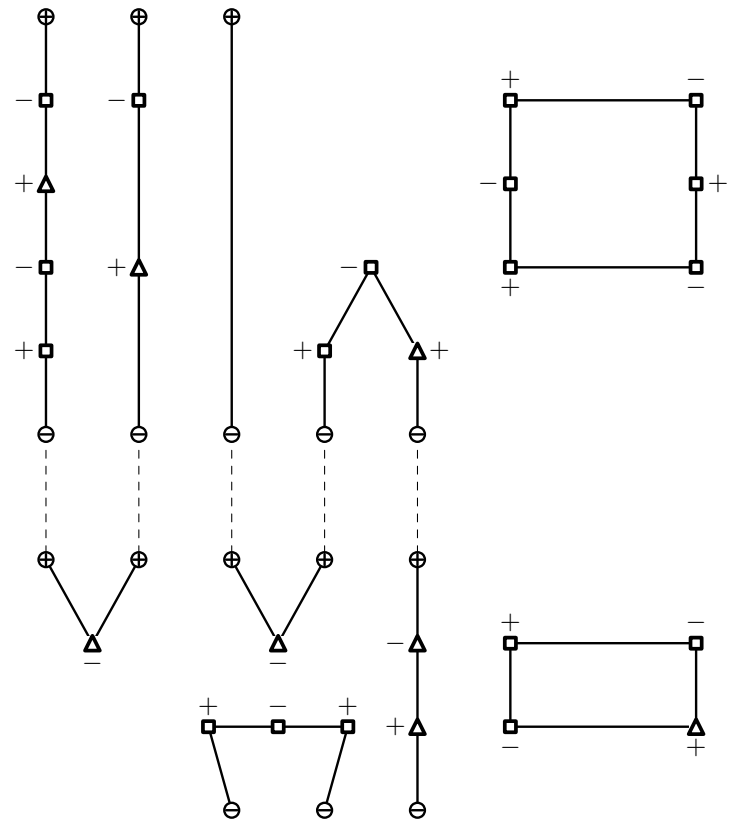}
$$

\caption{Reference to Subsection \ref{ss:chips-graph}.
Two chips  and their product.
 We connect exits of the upper diagram
with entries of  the lower diagram.%
\label{fig:product-chips}}
\end{figure}

To define the product
of morphisms
$\frG:\Omega_{[\alpha]}\to \Omega_{[\beta]}$,
$\frH:\Omega_{[\beta]}\to \Omega_{[\gamma]}$,
we identify exits of the  diagram $\frG$ with
 the corresponding entries of $\frH$.
Points of gluing become interior points of 
edges, see Figure \ref{fig:product-chips}.


\sm

\begin{figure}
 $$
\epsfbox{chip.1}
$$
a) An element $g$ of $\S_\infty(\Omega)$ for $p=2$.
 Elements of  subsets
$\Omega_{[\alpha]}'$ are marked by $\epsfbox{chip.9}$,
elements of $\Omega_{[\beta]}^{''}$ by $\epsfbox{chip.10}$.
Horizontal arcs mark the involution in $\Omega$.

$$
\epsfbox{chip.2}
$$

b) The diagram corresponding to the element $g$.
We contract  horizontal arcs to a vertex
 and remember  smells
of these arcs (symbols  $\epsfbox{hurwitz.7}$, 
$\epsfbox{hurwitz.8}$).
We also set a symbol ``$+$'' (respectively, ``$-$''),
 which bears in mind  was an arc upper or lower.

\caption{Reference to Section 4.\label{fig:3}}

\end{figure}

\sm

{\bf \punct Construction of a diagram from an element
of $\S_\infty$.%
\label{ss:chips-chips}}
For each $\omega\in\Omega$, we draw 
a segment with vertices $\omega'$ and $(g\omega)^{''}$. 

\sm

1. We identify   $\omega'$
with $(\omega^\circ)'$ and 
$\omega^{''}$ with $(g\omega^\circ)^{''}$.
Thus we get a countable disjoint union of finite chains.
Endpoints of open chains are elements of
$\Omega_{[\alpha]}'$ and $\Omega_{[\beta]}^{''}$.

\sm

2. We forget point $\omega$, $\omega^{''}$ corresponding to  
interior points of chains and remember only their smells.

\sm

3. We remove cycles of length 2 having
 vertices of the same smell.

\sm

See Figure \ref{fig:3}.


\sm


{\bf \punct Examples of representations.%
\label{ss:rep-chips}}
For a Hilbert space $H$ denote by
$$
{\mathsf S} H^k\subset H^{\otimes k}
$$
its symmetric power.

Consider a Hilbert space $V$,
 and collection of unit vectors
$$
\xi_i\in \mathsf S^2 V\subset V\otimes V
,$$ where $i=1$, \dots, $p$.
 Consider the tensor product
$$
\bigotimes_{i=1}^p
\left(\bigotimes_{\omega\in \N_i} (V\otimes V,\xi_i)
\right)\otimes \bigotimes_{\lambda\in L} V
$$
The group $G\cdot \bfK$ acts by permutations of factors $V$.

The parameter of a representation is a collection of vectors
$\xi_i\in \mathsf S^2 V$ defined up to action of 
the complete unitary group
$U(V)$ of the Hilbert space $V$.

Note that for $p=1$ (the case considered by Olshanski 
\cite{Olsh-symm}) we have a unique vector $\xi$ and it can be reduced
to a diagonal form.

\section{Chips (continuation)}

\COUNTERS

Here we discuss the same group $G$ and change the
subgroup. Also, we change language (pass from a graph 
to a dual graph).

\sm


{\bf\punct Group.%
\label{ss:chips-bis}} Let
 $Z=\begin{pmatrix}2&\dots &2 \end{pmatrix}$.
Consider the same objects
$\N_j$, $\Omega$, $L$, as in Subsection \ref{ss:chips1}.
Consider the pair
$$
(G,K)=(G[Z,\Lambda],K^\circledast[Z])
,$$
i.e., 
$$
K=\prod_{i=1}^p \S_\infty(\N_i)
$$

Now for each $\N_i$ we have two embeddings 
$\N_i\to\Omega$, which are distinguished by their melodies.

Consider the same sets $\Omega_{[\alpha]}$
as in  Subsection \ref{ss:chips1}.

\sm

{\it The  difference is  following.}
In the previous section ends of horizontal arcs in
 Fig.\ref{fig:3}.a  had equal rights.
Now they 
have different melodies.
Since there are only two melodies, we prefer to think that 
 arcs are oriented and arrows look from the left to right
(for instance, we think that   origins 
of arcs are 'violins' 
and  ends are 'contrabass's').
Now we contract {\it vertical} (i.e., inclined) arcs
and get a graph composed of chains and cycles.

We explain the construction more carefully.

\sm

\begin{figure}
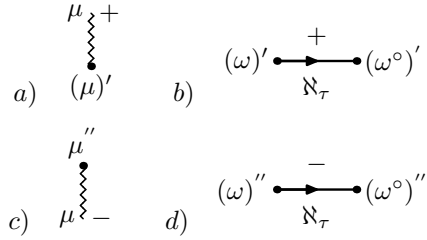

$$
a)\quad \epsfbox{chip.11}\qquad
b)\quad \epsfbox{chip.12}\qquad
$$
$$
c)\quad \epsfbox{chip.13}\qquad
d)\quad \epsfbox{chip.14}\qquad
$$
\caption{Reference to Subsection \ref{ss:chips-graph}. 
Items of a chip\label{fig:items}}
\end{figure}


{\bf \punct Construction of the graph.%
\label{ss:chips-graph}} Consider the following
 items for building 
a graph, see Figure \ref{fig:items}.

\sm

a) Let $\mu\in \Omega_{[\alpha]}$.
We draw a vertex and a tail as it is shown 
on Figure \ref{fig:items}. We write $\mu$ and ``$+$''
on the tail. Also  we write  a label
$\mu'$ on the vertex.

\sm

b) Let $\omega\notin \Omega_{[\alpha]}$. For definiteness,
let the melody of $\omega$ be  'violin', therefore
the melody of $\omega^\circ$ be 'contrabass'.
Then we draw an oriented segment with 
the label
$\omega'$ on the origin and the label $(\omega^\circ)^{'}$
on the end.
We  mark these segments by the smell of $\omega$ 
(it coincides with
the smell of $\omega^\circ$). We also write 
the label ``$+$'' on the segment.

\sm

c)  Let $\nu\in \Omega_{[\beta]}$.
We draw a vertex and a tail. We write $\mu$ and ``$-$''
on the tail. Also  we set  a label
$\nu^{''}$ on the vertex.

\sm

d) Let $\omega\notin \Omega_{[\beta]}$.
Let the melody of $\omega$ be  'violin'.
Then we draw an oriented segment with 
the label
$\omega^{''}$ on the origin and the label 
$(\omega^\circ)^{''}$ on the end.
We  mark these vertices by the smell of $\omega$.
 We also write 
the label ``$-$'' on the segment.


\sm

Thus we get a collection of items.
For each $\omega\in\Omega$ we identify the vertex
with label $(\omega)'$ and 
the vertex with label $(\omega)^{''}$.

Next, we forget labels on vertices.

\sm

\begin{figure}
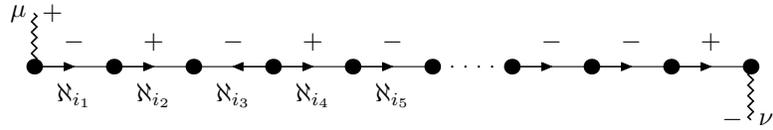
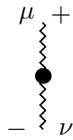
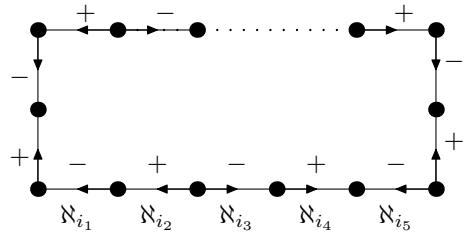
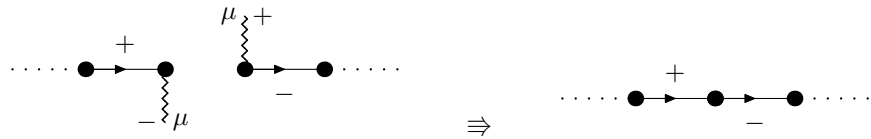

 $$ \epsfbox{chip.15}$$

a) An open chain.

 $$ \epsfbox{chip.17}$$

b) A short open chain.

 $$ \epsfbox{chip.16}$$

c) A closed chain.

$$\epsfbox{chip.18}\qquad \Rrightarrow\qquad  \epsfbox{chip.19}
$$
d) Gluing of chains.  

\caption {Reference to Section 6.\label{fig:chains-bis}}
\end{figure}

We get collection of chains of two types
(see Figure \ref{fig:chains-bis}).

\sm

a) {\it An open chain.} It is a chain of oriented segments,
each segment has a smell and a sign. An end of a chain is
a tail
(equipped with a label $\mu$
 and a sign). Signs  interlace.

We say that end vertices with ``$+$`` tails are {\it entries},
 and vertices with
''$-$`` tails are {\it exits}.

\sm

b) {\it A closed chain}
 consists of oriented segments equipped with
smells and signs. Signs  interlace.

It remains to define the multiplication
of two chips $\frG:\Omega_{[\alpha]}\to \Omega_{[\beta]}$,
$\frh:\Omega_{[\beta]}\to \Omega_{[\gamma]}$.
We identify exits of $\frG$ with corresponding  entries of $\frH$
and cut off tails.

\begin{proposition}
 This multiplication coincides with the multiplication in the train of
$(G,K)$.
\end{proposition}

It seems that this is self-obvious, see Figure 
\ref{fig:razdvig}. \hfill $\square$

\sm


{\bf\punct Examples of representations.%
\label{ss:rep-chip-bis}} We repeat 
the construction of Subsection \ref{ss:rep-chips}.
Now we can choose  distinguished vectors 
$\xi_j\in V\otimes V$ in arbitrary way 
(not necessary $\xi_j\in \mathsf S V^2$).

\section{Example: triangulated surfaces}

\COUNTERS

{\bf\punct Group.%
\label{ss:group-tri}}
Let $Z=\begin{pmatrix} 3\end{pmatrix}$,
$\lambda\ge 0$ be arbitrary. First, we consider
the pair
$$(G(Z,\Lambda), K^\circ(Z))=
\left(\S_{\lambda+3\infty},
 \S_\infty\ltimes (S_3)^\infty\right)
.$$
We reduce the subgroup and assume 
$$
K:=\S_\infty\ltimes (\Z_3)^\infty
,$$
where $\Z_3\subset S_3$ is the group the group of cyclic permutations 
(or, equivalently, the group of even permutations).

Now
$$
\Omega=L\sqcup \left(\N\times\I(3) \right)
.$$

Let $\alpha\ge 0$, the set $\Omega_{[\alpha]}$ 
is defined  as above (as fixed points of $K^\alpha$),
$$
\Omega_{[\alpha]}=L\sqcup (\M(\alpha)\times \I(3))
.$$

Thus we have only one color, only one smell, but 3 melodies,
say harp ($\nabla$), violin ($\heartsuit$), tube 
($\succ$).

\sm

\begin{figure}
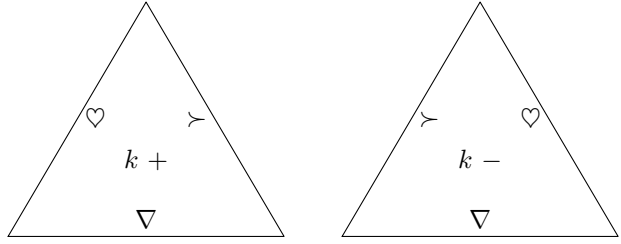
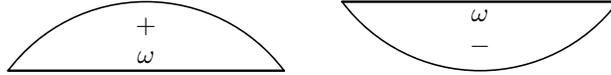


$$
\epsfbox{belyi.1}\qquad\epsfbox{belyi.2}
$$
a) A plus-triangle and a minus-triangle.

$$
\epsfbox{belyi.3}\qquad\epsfbox{belyi.4}
$$
b) A plus-tag and a minus-tag.

\caption{Reference to Section 7. 
 Items for a complex.\label{fig:items-treugolnik}
} 
\end{figure}

\begin{figure}
$$\epsfbox{belyi.5}$$

 A piece of complex. Removing numbers and melodies
$\succ$, $\nabla$, $\heartsuit$ (and leaving signs and labels
on tags), we
pass to double cosets.

$$\epsfbox{belyi.9}$$

\caption{Reference to Section 7.\label{fig:3complex}}
\end{figure}

\begin{figure}
 $$\epsfbox{belyi.6}$$

a) A degenerate component. An edge with two tags.

$b)\qquad\epsfbox{belyi.7}\qquad c)\qquad\epsfbox{belyi.8}$ 

 This a stereographic projection of a sphere and a graph on 
a sphere.  A pure  envelope (b) and an envelope (c).

\caption{Reference to Section 6. \label{fig:envelopes}}
\end{figure}


{\bf\punct Encoding of elements of symmetric
group.%
\label{ss:complex-tri}}
Fix $\alpha$, $\beta\ge 0$.

First, we take the following collection of items
(see Figure \ref{fig:items-treugolnik}).

\sm

A. {\it Plus-triangles and-minus triangles.}
 For each element of $k\in \N$
 we draw a pair of
oriented  triangles $T_\pm(k)$ with label $k$.  
We write labels $\nabla$, $\heartsuit$, $\succ$ 
on the sides $T_+(k)$ (resp. $T_-(k)$) clock-wise
(resp. anti-clock-wise). 

An important remark:
a number $k$
and a melody determines some element of $\Omega$.

\sm

B. {\it Plus-tags and minus tags.}
 For each element $\omega\in \Omega$ we draw two oriented
segments 
$D_\pm(\omega )$ with tags,
see Figure
\ref{fig:items-treugolnik}. We write the label $\omega$ and label '$+$'
(respectively '$-$')  on the segment 
$T_+(\omega)$ (resp. $T_-(\omega)$).

In this way, we have produced too much items.
Next, we  remove all  

\sm

--- triangles $T_+(k)$ with $k\le \alpha$;

\sm

--- triangles $T_-(k)$ with $k\le \beta$;

\sm

--- tags $D_+(\omega)$, where 
$\omega\notin\Omega_{[\alpha]}$;

\sm

--- tags $D_+(\omega)$, where 
$\omega\notin\Omega_{[\beta]}$.

\sm

Now each element of $\Omega$ is present
on precisely one edge of one item $T_+(k)$ or $D_+(\omega)$
(and respectively on one item $T_-(k)$ or $D_-(\omega)$).

Fix element $g\in \S_\infty(\Omega)$.
For each $\omega\in \Omega$ we identify
(keeping in mind the orientations) 
the edge of $T_+(\cdot)$ or $D_+(\cdot)$ labeled by 
$\omega$ with the edge of $T_-(\cdot)$ or $D_-(\cdot)$
labeled by $g\omega$.

In this way, we get a two-dimensional oriented triangulated
 surface $\Xi(g)$
with tags on the boundary. Our picture satisfies
the following properties: 

\sm

(i)  A surface consists of 
a countable number of compact  components.

\sm

(ii) Each component is a  two-dimensional oriented triangulated
 surface 
with tags on the boundary (we admit also
a segment with two tags, see Figure \ref{fig:envelopes}.a).

\sm

(iii) All triangles have labels '$+$' or '$-$',
neighboring triangles have different signs.

\sm

(iv)  Plus-triangles (resp. minis-triangles) are enumerated
by $\alpha+1$, $\alpha+2$, $\alpha+3$, \dots
(resp. $\beta+1$, $\beta+2$, $\beta+3$, \dots).

\sm

(v) Sides of plus-triangles are labeled (from
interior) by $\nabla$, $\heartsuit$, $\succ$ 
clockwise. Sides of minus-triangles are labeled by the same symbols
anti-clock-wise.

\sm

(vi) Tags are labeled by $\pm$. Plus-tags are enumerated
by elements of $\Omega_{[\alpha]}$, minus-tags
by elements of $\Omega_{[\beta]}$.

\sm

(vii) Almost all components are spheres 
composed of two triangles and labels on sides of the triangle 
coincide. 


\sm

We consider such surfaces up isotopes preserving
the orientation.

\begin{lemma}
Each surface equipped with data given above has the form
$\Xi[g]$. Different $g\in \S_\infty(\Omega)$
 produce different equipped surfaces.  
\end{lemma}

{\sc Proof.} We present the inverse construction. Above we 
have assigned two elements of $\Omega$ to each edge.
Let $\mu$ corresponds to the plus-side, $\nu$ corresponds
to the minus-side. Then $g$ send $\mu$ to $\nu$.
\hfill $\square$

\sm

Thus we get a bijection.

\sm

Now we need a technical definition.
We say that an {\it envelope} is a component consisting
of two triangles.
We say that an envelope is {\it pure} 
if melodies on both sides of each edge coincide 
(see Figure \ref{fig:envelopes}).

\sm


{\bf\punct Projection to double cosets.}

\begin{lemma}
\label{l:cleaning}
  Right multiplications $g\mapsto gh$ by elements of 
$\S^\alpha_\infty(\N)$ correspond to permutations of labels
$\alpha+1$, $\alpha+2$, $\alpha+3$, \dots on plus-triangles.
Respectively left multiplication by elements
of $\S^\beta_\infty(\N)$ correspond to permutations
of labels on minus-triangles.
\end{lemma}

\begin{lemma}
 Right multiplications by elements of the group
$(\Z_3)^{-\alpha+\infty} 
\subset 
K^\alpha$ correspond to cyclic permutations of symbols 
$\nabla$, $\heartsuit$, $\succ$ inside 
each plus-triangle.
\end{lemma}

\begin{corollary}
Pass to double cosets $K^\beta\setminus G/K^\alpha$
corresponds to forgetting numbers of triangles, 
melodies of sides, and removing all  envelops.

\end{corollary}

Thus, we remember only labels on tags and signs.


\sm

{\bf\punct Construction of the train.%
\label{ss:train-tri}}
Objects of the category are  $\alpha\ge 0$.
Fix indices $\alpha$ and $\beta$.
A {\it morphism} $\alpha\to\beta$
is a compact (generally, disconnected) triangulated 
surface without envelopes equipped with 
data (iii), (vi) from the list above
(labels $\pm$ and labels on tags).

To multiply $\frG:\alpha\to\beta$,
$\frH:\beta\to\gamma$, we glue (according
the orientations) minus-segments
of the boundary of $\frG$ with plus-segments
of the boundary of $\frH$ having the same labels.
We remove corresponding tags, forget their labels, 
forget the contour of gluing. It can appear some envelops,
 we remove them.

We get a morphism $\alpha\to\gamma$.

\begin{theorem}
\label{th:tri-train}
The multiplication described above is the multiplication
in the train of the pair
$(G,K)$. 
\end{theorem}

{\sc Proof.} Fix $h$, $g\in\ G$.
Consider the corresponding
$\Xi[h]$, $\Xi[g]$. Take very large 
$N$. Then the set of labels $k$ on minus-triangles
of $\Xi[\theta_\beta^N h]$ and the set of labels $l$ on
plus-triangles of $\Xi[g]$ are disjoint.
Therefore plus-triangles of $\Xi[\theta_\beta^N h]$  
preserve their neighbors after the multiplication
 $\Xi[\theta_\beta^N h]\to g\Xi[\theta_\beta^N h]$.
 Also, minus-triangles of $\Xi[g]$ preserve their neighbors
after multiplication $g\mapsto  g\Xi[\theta_\beta^N h]$.
Therefore both surfaces $\Xi[h]$, $\Xi[g]$
are pieces of the surface $\Xi[g\theta_\beta^N h]$.
\hfill $\square$

\sm


{\bf\punct Involution on the train.}
We  reverse signs and reverse the orientation.

\sm


{\bf\punct Examples of representations.%
\label{ss:reps-tri}}
Let $V$ be a Hilbert space.
Fix unit vectors
\begin{equation} 
\xi_1,\, \xi_2,\, \xi_3\in V\otimes V\otimes V
\label{eq:s3}
\end{equation}
invariant with respect to cyclic permutations
 of elements of the tensor products.

Consider the tensor product
$$
\bigotimes_{l\in L} V 
\otimes
\bigotimes_{i=1}^3  (V\otimes V\otimes V,\xi_i)^{\otimes\infty}
.$$
The group
$G\cdot\bfK$ acts on this product by permutations of factors.

\sm


{\bf\punct Another pair.%
\label{ss:tri-bis}} Let $Z$, $\lambda$, $G=G[Z,\lambda)$
 be the same,
consider the pair
$$
(G,K):=(G(Z,\lambda), K^\circledast(Z))=
(\S_{\lambda+3\infty}, \S_\infty)
.$$
Return to Lemma \ref{l:cleaning}.
Now we remove numbers of triangles but preserve
melodies. We also remove all pure envelopes.

In construction of representation, we can replace
(\ref{eq:s3}) by arbitrary
unit vectors
$$
\xi_1,\, \xi_2,\, \xi_3\in V\otimes V\otimes V.
$$


{\bf \punct  Another pair.}
Let $Z$, $\lambda$, $G=G[Z,\lambda)$
 be the same.
Consider the pair
$$
(G,K):=(G(Z,\lambda), K^\circ(Z))=
(\S_{\lambda+3\infty}, \S_\infty\ltimes (S_3)^\infty)
.$$
First, we construct representations. In the construction
of Subsection \ref{ss:reps-tri} we take 
$$
\xi_1,\, \xi_2,\, \xi_3\in \mathsf S^3 V
\subset V\otimes V\otimes V.
$$

An attempt to repeat the  construction of the train
meets an obvious difficulty: permutations of melodies
change orientations of triangles. However, 
we can pass from triangulations to dual graphs.
Now we can enumerate double cosets by tri-valent graphs.
See the following section.


\sm

\section{General case, $K=K^\circ$ is a wreath product}

\begin{figure}
$$ \epsfbox{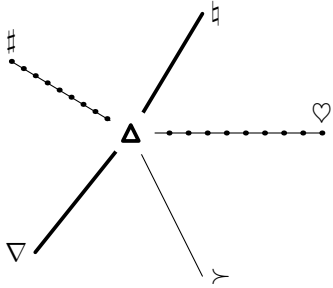} $$

a) A node. 

$$ \epsfbox{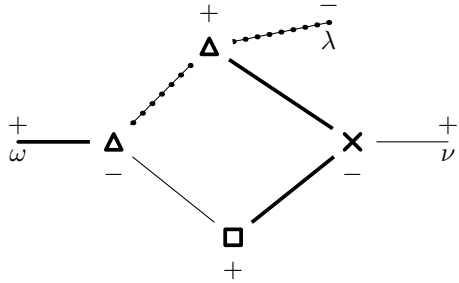} 
$$

b) A double coset.

$$ \epsfbox{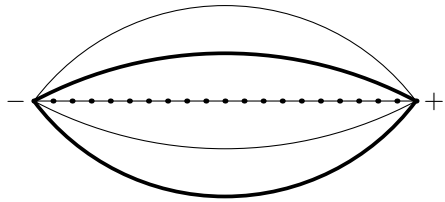} 
$$

c) A {\it trivial} component of a graph.

\caption{Reference to Subsection \ref{ss:coding-2}.
\label{fig:graph}}
\end{figure}

\COUNTERS

{\bf\punct Group.%
\label{ss:general-3}} Here we consider an arbitrary matrix
$Z$ and an arbitrary vector $\Lambda$.
Now
\begin{align*}
\Omega_j&=L_j\sqcup \coprod_{i=1}^{p}
\bigl(\N_i\times \I(\zeta_{ji})\bigr),
\\
G&:=G[Z,\Lambda]=\prod_{j=1}^q \S_\infty(\Omega_j),
\\
K&:=K^\circ[Z]=\prod_{i=1}^p \Bigl(\S_\infty(\N_i)\ltimes 
\Bigl(\prod_j S_{\zeta_{ji}}\Bigr)^\infty\Bigr)
\end{align*}
Recall that we have attributed a color to each $\Omega_j$,
a smell to each $\N_i$, and a melody
to each infinite orbit of $\S_\infty(\N_i)$ on
$\Omega_j$.

We denote 
$$
\Omega:=\coprod_{j\le q} \Omega_i
$$
and regard $G[Z,\lambda]$ as as a subgroup in $\S_\infty(\Omega)$.
For a multiindex $\alpha=(\alpha_1,\dots,\alpha_p)$
 we denote by 
$\Omega_{[\alpha]}$ the set of all 
$K^\alpha$-fixed points of $\Omega$,
$$
\Omega_{[\alpha]}=
\coprod_{j=1}^q \left(L_j\sqcup \coprod_{i=1}^p 
\bigl(\M(\alpha_i)\times \I(\zeta_{ji})\bigr)\right)
.$$


{\bf\punct Encoding of elements of  symmetric groups.%
\label{ss:coding-2}}
For each element of $G[Z,\lambda]$ we construct  
a graph equipped with some additional data.

\sm

 For each smell $i$ we draw a {\it node}
$T[\aleph_i]$ (see Figure {fig:graph}). It contains
a vertex of smell $\aleph_i$ and 
$\sum_j \zeta_{ji}$ semi-edges. Edges are colored,
each color $\gimel_j$ is used for coloring $\zeta_{ji}$ edges.
Also we attribute a melody to each semi-edge.
Thus semi-edges of $T[\aleph_i]$
 are in one-to-one correspondence with
orbits of $\S_\infty(\N_i)$ on $\Omega$.

Now we  prepare the following collection of items.

\sm

a) For each smell $i$ and each $k\in \N$ we draw
two copies  $T_\pm[\aleph_i;k]$ of the node $T[\aleph_i]$,
 their vertices
are  labeled $k$ and $\pm$.
We through out nodes $T_+[\aleph_i;k]$ with $k\le\alpha_i$
and $T_+[\aleph_i;m]$ with $m\le\beta_i$

\sm

b) For each color $j$ and for each element $\omega$ of 
$\Omega_j\cap \Omega_{[\alpha]}$ we draw
a tag $D_+(\omega)$ and mark this tag by $\omega$, the color of 
$\omega$, and ''$+$``. We draw similar tags 
$D_-(\omega)$ for elements  $\omega\in\Omega_{[\beta]}$.
We imagine a tag as a vertex and semi-edge.

\sm

Thus the set  $\Omega$ is in one-to-one correspondence 
with the sets 
$$
\cE_+=
\left\{ \begin{matrix}\text{ All semi-edges of}\\
\text{  all nodes $T_+[\aleph_i,k]$}
\end{matrix}
\right\}
\bigcup \Omega_{[\alpha]}
$$
and 
$$
\cE_-=
\left\{ \begin{matrix}\text{ All semi-edges of}\\
 \text{ all nodes $T_-[\aleph_i,k]$}
 \end{matrix}
\right\}
\bigcup \Omega_{[\beta]}
.$$
Denote the bijections $\Omega\to \cE_\pm$ by
$H_\pm$.

Now for each $\omega\in\Omega$ we connect a semi-edge  
$H_+(\omega)\in \cE_+$ with the semi-edge
$H_-(g\omega)\in \cE_-$. We get a graph
with following properties.

\sm

(i) Graph consists of a countable number of compact components.

\sm

(ii) There are 2 types of vertices, interior vertices%
\footnote{The case $\sum_j \zeta_{ji}=1$ is admissible,
then we can meet a unique edge in an interior vertex
 of the smell $\aleph_i$}
and end-vertices (ends of semi-edges).

\sm

(iii) Each interior vertex has a smell $\aleph_i$
and 
a sign ''$+$`` or ''$-$``.

\sm

(iv) Interior plus-vertices are enumerated by the set 
$\{\alpha+1,\alpha+2,\dots\}$, interior minus-vertices by
$\{\beta+1,\beta+2,\dots\}$.

\sm

(v) End-vertices fall into two classes, entries and exits.
Entries are enumerated by  elements of $\Omega_{[\alpha]}$
and labels ''$+$``. Exits are enumerated by 
 elements of $\Omega_{[\beta]}$
and labels ''$-$``

\sm

(v)  Neighboring vertices have different signs.

\sm

(vi)  Edges are colored, number of edges of a color $\gimel_j$
coming to an interior vertex of smell $\aleph_i$ is $\zeta_{ji}$.
The edge adjacent to an end vertex with label $\omega$
 has the color of $\omega$.

\sm

(vii) For each semi-edge adjacent to an interior
vertex it is attributed a melody 
compatible with its
color. At each vertex of smell $\aleph_i$ 
each compatible with the smell melody is present precisely
one time.

\sm

(viii) All but finite number of components 
consists of two interior vertices and edges connected this vertices.

\sm

We call components described in (viii) {\it trivial}. 
We say that a component is {\it completely trivial}
if for each edge smells of both semi-edges coincide.

\begin{theorem}
There is one-to-one correspondence between
the set of all graphs satisfying {\rm (i)-(viii)} 
and the infinite symmetric group 
\end{theorem}

{\sc Proof.} Consider an edge. It has a plus-semi-edge
and a minus-semi-edge. Consider the corresponding elements
$\phi\in\cE_+$ and $\psi\in\cE_-$.
We set $g\phi=\psi$.
\hfill $\square$

\sm


{\bf\punct Projection to double cosets.}

\begin{proposition}
{\rm a)} Right multiplications by elements of 
$\S^\alpha_\infty(\N_i)$ correspond to permutations of labels 
$\{\alpha+1,\alpha+2,\dots\}$ on plus-vertices
of the smell $\aleph_i$.

\sm

{\rm b)} Right multiplications by elements
of $S_{\zeta_{ij}}^{-\alpha_i+\infty}\subset \S_\infty(\Omega_j)$
correspond to permutations of melodies of
 semi-edges of color $\beth_j$
adjacent to fixed vertices of the smell
$\aleph_i$.
 \end{proposition}

\begin{corollary}
Projection to double cosets correspond to forgetting labels
$\in\N$ and melodies 
\end{corollary}

Colors, smells, signs, and also labels on tags are preserved.

\sm


{\bf\punct Multiplication of double cosets.}
For two morphisms $\frG:\alpha\to\beta$,
$\frH:\beta\to \gamma$, we glue exits of $\frg$
 with corresponding entries of $\frh$.

The involution is the inversion of signs and also entries/exits.

\begin{theorem}
 This product coincides with the product
in the train $\T(G,K)$.
\end{theorem}

{\sc Proof.} See proof of Theorem \ref{th:tri-train}.

\sm

{
\bf\punct Some representations of $(G,K)$.%
\label{ss:rep-general}}
We consider a collection of Hilbert spaces $W_1$,
\dots, $W_q$ enumerated by colors.
Fix $i$. Consider the tensor product
\begin{equation}
\cH_i=
\bigotimes_{j=1}^p W_j^{\otimes \zeta_{ji}}
\label{eq:cH}
.\end{equation}

\begin{observation}
\label{obs}
Factors of the product are in one-to-one correspondence
with semi-edges a
$T[\aleph_i]$
\end{observation}

Fix a unit vector  
\begin{equation}
\xi_i\in\bigotimes_{j=1}^p 
\mathsf S^{\zeta_{ji}}W_j\subset \cH_i
\label{eq:111}
.
\end{equation}

Consider the tensor product
$$
\frW:=\bigotimes_{j=1}^q W_j^{\otimes \lambda_j}
\otimes \bigotimes_{i=1}^q (\cH_i,\xi_i)^{\otimes \infty}
.$$
Note that factors $W$ of this tensor product are enumerated by
elements of $\sqcup \Omega_j$.
Formally we can write
$$
\bigotimes_{j=1}^q \bigotimes_{\omega\in\Omega_j} W_j
$$
However, this makes no sense without distinguished vectors.

Each group $\S_\infty(\Omega_j)$ acts by permutations of factors
$W_j$. This determines the action of $G$.

The groups $\S_\infty(\N_i)$
act as permutations of factors in
$$(\cH_i,\xi_i)^{\otimes \infty}.$$

Thus we get action of $G\cdot\bfK$.
 
\sm


{\bf \punct Reduction of the subgroup.} 
Now consider the pair $(G,K)=(G[Z,\Lambda],K^\circledast[Z])$.
We use the same encoding of $G[Z,\Lambda]$.

For passing to double cosets 
$K^\beta\setminus G/K^\alpha$ we forget numeric labels 
(but remember melodies). Also we remove completely trivial
components of the graph.

In the following section we give another (equivalent)
description of double cosets.

\section{General case, $K$ is a product of symmetric groups}

\COUNTERS

Construction of this section is more-or less a version
 of the previous construction.
For smaller group $K$ we can replace a graph 
by a fat graph and after this  draw a two dimensional surface.
We repeat the construction independently.

\begin{figure}
$$\epsfbox{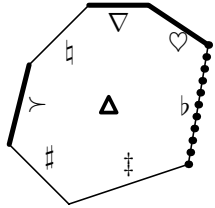} $$

\caption{Reference to Section 7. A polygon $T_+[\aleph_i]$.
\label{fig:polygon}}
\end{figure}

\begin{figure}
$$\epsfbox{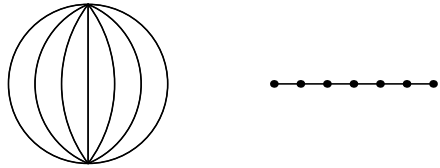}$$

\caption{Reference to Subsection \ref{ss:simple}. A di-gonal complex and
the corresponding one-dimensional chain.
\label{fig:digon}}
\end{figure}

\begin{figure}
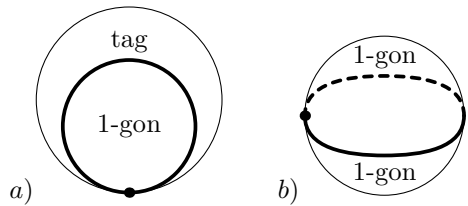

$$a) \epsfbox{last.3}\qquad b)  \epsfbox{last.4}$$

\caption{Reference to Section \ref{ss:simple}. The only
possible connected 1-gonal complexes.
\label{fig:onegon}}
\end{figure}

{\bf\punct Group.%
\label{ss:general-2}} We consider an arbitrary matrix
$Z$ and an arbitrary vector $\Lambda$.
Now
\begin{align*}
\Omega_j&=L_j\sqcup \coprod_{i=1}^{p}
\bigl(\N_i\times \I(\zeta_{ji})\bigr),
\\
G&:=G[Z,\Lambda]=\prod_{j=1}^q \S_\infty(\Omega_j),
\\
K&:=K^\circledast[Z]=\prod_{i=1}^p  \S_\infty(\N_i)
.
\end{align*}
Recall that we have attributed a color to each $\Omega_j$,
a smell to each $\N_i$, and a melody
to each infinite orbit of $\S_\infty(\N_i)$ on
$\Omega_j$.

We denote 
$$
\Omega:=\coprod_{j\le q} \Omega_i
.$$
For a multiindex $\alpha=(\alpha_1,\dots,\alpha_p)$
 we denote by 
$\Omega_{[\alpha]}$ the set of all 
$K^\alpha$-fixed points of $\Omega$,
$$
\Omega_{[\alpha]}=
\coprod_{j=1}^q \left(L_j\sqcup \coprod_{i=1}^p 
\bigl(\M(\alpha_i)\times \I(\zeta_{ji})\bigr)\right)
.$$


\sm

{\bf \punct Constructions of the train.%
\label{ss:train-general}} 
Fix a smell $\aleph_i$. Nontrivial orbits of
$\S_\infty(\N_i)$ on $\sqcup \Omega_j$ are enumerated by a
pair (color, melody). Total number of such orbits is
$$
\sum_j \zeta_{ij}
.$$
{\it We choose an  arbitrary cyclic order
on the set of such pairs} (the construction below 
depends on this choice). Next, we draw a polygon
$T_+[\aleph_i]$
of smell $\aleph_i$, whose sides are marked by
pairs (color, melody)
according the cyclic order. 
We also define  the polygon $T_-[\aleph_i]$,
whose sides are marked according the reversed cyclic order.

\sm

{\sc Remark.} Let us forget melodies on the sides
 of the polygon. If there is no rotation of the polygon
preserving the coloring of sides, then melodies
 can be reconstructed
 in a unique way. 
\hfill$\square$ 

\sm

Consider the following collection of items.

\sm

--- {\it Plus-polygons and minus-polygons.} For each $i\le p$
for each 
$k\in\N_i$ we draw the pair of oriented polygons $T_\pm[\aleph_i,k]$
that were described above;
they are additionally labeled by $k\in \N$.
Any side of any polygon
 $T_\pm[\aleph_i,k]$ has smell, color, and melody;
 therefore a side determines an element
of  $\Omega$.

\sm

--- {\it Plus-tags and minus-tags.}
For each element  $\omega\in\Omega$, we draw two tags
$D_\pm(\omega)$ labeled by $\omega$ and $\pm$, see Figure
\ref{fig:items-treugolnik}. The side of a tag is 
painted in the color of $\omega$

\sm

Fix multi-indices $\alpha$, $\beta$.

We remove some items from the collection:

\sm

--- polygons $T_+[\aleph_i,k]$ if $k\le \alpha_i$;

\sm

--- polygons $T_-[\aleph_i,m]$ if $m\le \beta_i$;

\sm

--- tags $T_+[\omega]$ if $\omega\notin \Omega_+{[\alpha]}$;

\sm

--- tags $T_-[\omega]$ if $\omega\notin \Omega_-{[\beta]}$.

\sm

Now we have the one-to-one correspondences between 
the set $\Omega$ and set of all edges of all plus-triangles
and plus-tags. Also we have  one-to-one correspondences between 
the set $\Omega$ and set of all edges of all minus-triangles
and minus-tags.

\sm

For each $g\in G$ we glue the complex. For each
 $\omega\in\Omega$, we identify
the (oriented) edge of a plus-polygon or a plus-tag corresponding
$\omega$ with the (oriented) edge of 
a minus-polygon or minus-tag corresponding $g\omega$.

\sm

Thus we get a polygonal two-dimensional oriented surface 
with tags on the boundary satisfying the following
properties:

\sm

(i) The surface consists of countable number of compact
components.

\sm

(ii) Each component is  tilled by polygons 
of the types $T_\pm[\aleph_i]$
and has tags $D_\pm$ on the boundary.

\sm

(iii) Each polygon is labeled by '$+$' or '$-$' ,
neighboring polygons have different signs.

\sm

(iv) Each edge has a color, which is common for both
(plus and minus sides of an edge).

\sm

(v) Each edge has two melodies, on the plus-side
and on the minus-side.

\sm

(vi) Cyclic order of pairs (color, melody) around
the perimeter of each polygon $T_\pm[\aleph_i]$ is fixed.

\sm

(vii) Plus-polygons (resp., minus-polygons) of a fixed smell
$\aleph_i$
are enumerated by $\alpha+1$, $\alpha+2$. \dots (resp.
$\beta+1$, $\beta_2$, \dots).

\sm

(viii) Plus-tags are enumerated by points of $\Omega_{[\alpha]}$
and minus-tags by points of $\Omega_{[\beta]}$.

\sm

(ix) All but a finite number of  components of the surface
 are unions of  pairs $T_+[\aleph_j,k]$ and $T_-[\aleph_j,l]$.
We call such components '{\it envelopes}'.
We say that a {\it pure envelope} is an envelope
such that melodies on plus and minus sides of each edge coincide.

\sm

\begin{theorem}
Data of such type are in one-to-one correspondence with
the group $G$.
\end{theorem}

{\it Inverse construction.}
 To find $g\omega$, we find $\omega$
inside pairs (plus-polygon, side). 
This side also is a side of minus-polygon and codes
the element $g\omega$.


\sm

{\bf\punct Pass to double cosets $K^\beta \setminus G/ K^\alpha$.}
The literal  analog of Lemma \ref{l:cleaning} holds.

To pass to double cosets
$K^\alpha\setminus G/K^\beta$
we forget labels $k\in \N$
and remove all envelops.

We get a compact surface tiled by polygons;

\sm

--- Polygons are equipped with signs $\pm$ and smells, 

\sm

--- Each edge is equipped with a color and a
pair of melodies on the negative side and 
the positive side of the edge
(coloring and melodization of edges of each polygon is fixed up
 to a cyclic permutation
of sides as above)%
\footnote{Recall that in many cases melodies can be uniquely 
reconstructed from colors and may be forgotten}. 

\sm

--- The boundary edges of the surface
are equipped with signs, positive edges 
 are enumerated by points of $\Omega_{[\alpha]}$,
negative edges  by points of $\Omega_{[\beta]}$.

\sm

We say that such surface is a morphism
$\alpha\to\beta$.

Let $\frG:\alpha\to\beta$, $\frH:\beta\to\gamma$
be two surfaces. For each $\omega\in \Omega{[\beta]}$
we glue $\omega$-exit of $\frG$
with $\omega$-entry of $\frH$ (according the orientation).
Removing envelops, we come to a complex of the same type.

\begin{theorem}
This multiplication coincides with the multiplication 
in the train, 
\end{theorem}

{\sc Proof} is the same as for Theorem \ref{th:tri-train}.


\sm

{\bf\punct Involution in the  train.}
We reverse sign $\pm$ and reverse the orientation.

\sm


\sm

{\bf \punct Simple cases.%
\label{ss:simple}} Note that our construction
admits 2-gons and 1-gons.
Let the matrix $Z$ satisfies 
$\sum_j \zeta_{ij}=2$ for all $i$. Then all our polygons
are 2-gons. A digonal complex can be regarded as union of chains
and we come to chip-type constructions (see Sections 4, 6).

If $Z=\begin{pmatrix}1&\dots &1 \end{pmatrix}$,
then our complex consists of $1$-gons. There are only
3 types of possible components and we come to construction of 
Section 3.

\sm


{\bf \punct Belyi data,%
\label{ss:belyi}} see \cite{Bel1}, \cite{Bel2},
\cite{LZ}, \cite{Wal}. Consider the Riemann  sphere
$\ov \C=\C\cup\infty$, compact closed Riemannian surface $\Gamma$ and
a covering map $\Gamma\to \ov \C$ having
ramifications only at points $0$, $1$, $\infty\in \ov\C$.
We draw intervals $[0,1]$, $[1,\infty]$,
$[\infty,0]$ on the sphere and paint them in red, yellow, blue
respectively. We set the label '$+$' on the upper 
half-plane in $\ov \C$
 and '$-$'
on the lower half-plane. 

Lift this picture to the covering $\Gamma$. We get a triangulation
of $\Gamma$, triangles are labeled by $\pm$ and edges are colored.
In our notation this corresponds
to double cosets
$$
K^\circledast[Z]\setminus G[Z,\Lambda]/K^\circledast[Z]
\qquad\text{
with $Z=\begin{pmatrix}
       1\\1\\1
      \end{pmatrix}
$, $\Lambda=0$}
$$

Belyi Theorem claims that such maps $\Gamma\to\ov \C$ exists if
and only if a curve $\Gamma$ is determined 
over algebraic closure of $\Q$.

\sm


{\bf \punct Constructions of representations.%
\label{ss:rep-general-bis}}
In the construction of tensor product
from Subsection \ref{ss:rep-general} we can choose  arbitrary
unit vectors
$$
\xi_i\in\bigotimes_{j=1}^p 
 W_j^{\otimes\zeta_{ji}}=: \cH_i.
$$
instead of (\ref{eq:111}).

\section{Spherical functions}

\COUNTERS

We wish to write explicit formulas for
spherical functions of some representations
in terms of trains. To be definite, we consider
pairs
$(G,K)=(G[Z,0],K^\circledast[Z])$
and a representation $\rho$ of $G$ in the space 
$$
\frW:=
 \bigotimes_{i=1}^q (\cH_i,\xi_i)^{\otimes \infty}
$$
described above.
The formula written below is a precise copy
 of \cite{Ner-symm}. We omit a proof.


\sm

{\bf\punct Formula.%
\label{ss:matrix-elements}}
We consider a morphism $\frG:0\to 0$,
i.e. a closed 
polygonal complex $\frG:\alpha\to\beta$.

Choose an orthonormal  basis $e_k[j]$ in each $W_j$.

Now we assign  elements $e_k[j]$ 
to edges of the complex
according
the following condition: 
colors of $e_k[j]$ correspond to colors of edges.

\sm

Look to  a polygon. Keeping in the mind  Observation \ref{obs},
we observe that a tensor product of basis vectors $e$ along
the perimeter of the polygon $P$ is a basis vector,
in $\cH$, say $E_\Phi(P)$.
The formula for spherical function is
\begin{equation}
\sum_\Phi
\prod\limits_{\text{plus-polygons $P$}}
\la E_\Phi(P), \xi_{\aleph(P)}\ra_\cH
\cdot
\prod\limits_{\text{minus-polygons $Q$}}
\la \xi_{\aleph(Q)}, E_\Phi(Q)\ra_\cH
\label{eq:matrix-element}
\end{equation}
where $\aleph(P)$ is the smell of $P$ and
the summation is given over all arrangements
of basis vectors.

\sm


{\sc Remark.} For chips this formula is reduced
to calculations of traces of products of matrices.


{\tt Math.Dept., University of Vienna,

 Nordbergstrasse, 15,
Vienna, Austria

\&

Institute for Theoretical and Experimental Physics,

Bolshaya Cheremushkinskaya, 25, Moscow 117259,
Russia

\&

Mech.Math.Dept., Moscow State University,

Vorob'evy Gory, Moscow

e-mail: neretin(at) mccme.ru

URL:www.mat.univie.ac.at/$\sim$neretin

wwwth.itep.ru/$\sim$neretin
}


\begin{thebibliography}{cc}

\bibitem{Ati}
Atiyah, M. {\it Topological quantum field theories.}
  Inst. Hautes \`Etudes Sci. Publ. Math.  No. 68  (1988), 175--186.

\bibitem{Bae}
 Baez, J. C. {\it Spin foam models.} 
 Classical Quantum Gravity  15  (1998),  no. 7, 1827--1858

\bibitem{Bel1}
Bely\u i, G. V.
 {\it Galois extensions of a maximal cyclotomic field.}
 Math. USSR Izv. 14 (1980), 247--256.


\bibitem{Bel2}
Bely\u i, G. V.
{\it A new proof of the three-point theorem.}
 Sb. Math.  193  (2002),  no. 3-4, 329--332.

\bibitem{Bra}
Brauer, R.
{\it On algebras which are connected with
 the semisimple continuous groups.} 
 Ann. of Math. (2)  38  (1937),  no. 4, 857--872. 

\bibitem{Dix}  Dixmier, J.
{\it Les $C^*$-algebr\`es et leurs repr\'esentations.}
Gauthier Villars, 1964.

\bibitem{Gui}
Guichardet, A.
{\it Symmetric Hilbert spaces and related topics. 
Infinitely divisible positive definite functions. 
Continuous products and tensor products.
 Gaussian and Poissonian stochastic processes.}
 Lecture Notes in Mathematics,
 Vol. 261. Springer-Verlag, Berlin-New York, 1972. 

\bibitem{Hur}
Hurwitz, A.,
{\it Ueber Riemann'sche Fl\"achen mit gegebenen Verzweigungspunkten.}
Math. Ann. 39 (1891), no. 1, 1--60. 

\bibitem{Ism1}
Ismagilov, R.S.,
{\it Elementary spherical functions on the groups $\SL(2,P)$ over
a field $P$, which is not locally compsct with respect to
the subgroup of matrices with integral elements.} 
Math. USSR-Izvestiya, 1967, 1:2, 349--380. 


\bibitem{Ism2}
 Ismagilov, R.S., {\it Spherical functions over 
a normed field whose residue 
field is infinite.} 
  Funct. Anal. and Its Appl.
Volume 4, N. 1, 37--45.


\bibitem{KOV}
Kerov, S., Olshanski, G., Vershik, A.
{\it Harmonic analysis on the infinite symmetric group.}
Invent. Math. 158 (2004), No. 3, 551--642.

\bibitem{Ker}
Kerov, S. V. {\it Realizations of representations
 of the Brauer semigroup.}
 Zap. Nauchn. Sem. Leningrad. Otdel. Mat. Inst. Steklov. (LOMI)  164
  (1987), 
 188--193; 
 translation in  J. Soviet Math.  47  (1989),  no. 2, 2503--2507

\bibitem{Kir}
 Kirillov, A. A. {\it Elements of the theory of representations.}
 Translated from the Russian by Edwin Hewitt.
Springer-Verlag, Berlin-New York, 1976.

\bibitem{LZ}
Lando, S. K., Zvonkin, A. K. 
{\it Graphs on surfaces and their applications.}
 Encyclopaedia of Mathematical Sciences,
 141. Low-Dimensional Topology, II. Springer-Verlag, Berlin, 2004

\bibitem{Lie}
Lieberman, A.
{\it The structure of certain unitary representations of
 infinite symmetric groups.}
Trans. Am. Math. Soc. 164, 189-198 (1972).


\bibitem{Nat}
Natanzon, S. M.
{\it Cyclic Foam Topological Field Theories.} 
 J. Geom. Phys.  60  (2010),  no. 6-8, 874--883. 


\bibitem{Ner-sto}
Neretin, Yu. A.
{\it 
Categories of bistochastic measures and representations of some infinite-dimensional groups.}
Russian Acad. Sci. Sb. Math. 75 (1993), no. 1, 197--219. 

\bibitem{Ner-book}
   Neretin, Yu. A. {\it Categories of symmetries and
 infinite-dimensional groups.}
 Oxford University Press, New York, 1996.

\bibitem{Ner-symm}
Neretin, Yu. A., {\it Infinite tri-symmetric group, multiplication 
of double cosets, and checker topological field theories.}
Preprint, available via {\tt arXiv:0909.4739}.

\bibitem{Ner-char} 
Neretin, Yu. A. {\it Multi-operator colligations and multivariate spherical functions.}
Preprint, available via {\tt arXive:1006.2275}.

\bibitem{Ner-faa}
Neretin, Yu. A. {\it
Sphericity and multiplication of double cosets for infinite-dimensional classical groups.} Preprint, available via {\tt arXiv:1101.4759} 

\bibitem{Ness1}
 Nessonov, N. I.   {\it Factor-representation of the group $\GL(\infty)$
and admissible representations $\GL(\infty)^X$} (Russian)
  Math. Phys., Analysis,
Geometry, 2003, N4,  167--187.

\bibitem{Ness2}
Nessonov, N. I.
Factor-representation of the group $\text{$\GL$}(\infty)$
and admissible representations of $\text{$\GL$}(\infty)^X$. II.
(Russian)
Mat. Fiz. Anal. Geom. 10, No. 4, 524--556 (2003).


\bibitem{Olsh-tree}
Olshanski, G. I.
{\it New "large" groups of type I.}
(Russian)
Itogi Nauki Tekh., Ser. Sovrem. Probl. Mat. 16, 31-52 (1980).
English transl.: 
J. Sov. Math. 18, 22--39 (1982).

\bibitem{Oko} Okounkov, A.
{\it Thoma's theorem and representations of the infinite
 bisymmetric group.}
 Funct. Anal. Appl. 28, No.2, 100--107 (1994).

\bibitem{Olsh-lieb}
Olshanski, G. I.
{\it Unitary representations of the infinite symmetric group:
a semigroup approach} in
{\it Representations of Lie groups and Lie algebras,}
Proc. Summer Sch., Budapest 1971, Pt. 2, 181--197 (1985).

\bibitem{Olsh-Opq}
Olshanski, G.I
{\it Infinite-dimensional groups of finite $R$-rank:
 description of representations 
and asymptotic theory.}
 Funct. Anal. and Appl., 1984, 18:1, 22--34.

\bibitem{Olsh-dop}
Olshanski, G. I.
{\it Unitary representations of the group ${\rm SO}_0(\infty,\infty)$
 as limits of unitary representations of the groups ${\rm SO}_0(n,\infty)$ as $n\to \infty$}. Funct. Anal. Appl. 20 (1986), no. 4, 292--301. 

\bibitem{Olsh-add1}
Olshanski, G. I. 
{\it The method of holomorphic extensions in the theory of unitary representations of infinite-dimensional classical groups.} 
  Funct. Anal. Appl.  22  (1988),  no. 4, 273--285 (1989)

\bibitem{Olsh-GB}
Olshanski, G.I.
{\it Unitary representations of infinite dimensional pairs
$(G,K)$ and the formalism of R. Howe.}
In
{\it Representation of Lie groups and related topics,}
Adv. Stud. Contemp. Math. 7, 269--463 (1990).

\bibitem{Olsh-symm}
 Olshanski, G.I.,
{\it Unitary representations of $(G,K)$-pairs connected with
the infinite symmetric group $S(\infty)$.}
Leningr. Math. J. 1, No.4, 983--1014 (1990).

\bibitem{Olsh-topics}
Olshanski, G. I. {\it On semigroups related
 to infinite-dimensional groups.}
In: {\it Topics in representation theory} (A. A. Kirillov,
ed.). Advances in Soviet Math., vol. 2. Amer. Math. Soc.,
Providence, R.I., 1991, 67-101.


\bibitem{vN}
von Neumann, J.
{\it On infinite direct products.}
 Compos. Math. 6, 1-77 (1938).
Reprinted in von Neumann {\it Collected works}.

\bibitem{Seg}
Segal, G. B.
 {\it  The definition of conformal field theory} In
{\it  Differential geometrical methods in theoretical physics 
(Como, 1987)}, 
 165--171, NATO Adv. Sci. Inst. Ser. C Math. Phys. Sci., 250,
 Kluwer Acad. Publ., Dordrecht, 1988.

\bibitem{Tho}
Thoma, E.
{\it Die unzerlegbaren, positiv-definiten
Klassenfunktionen der abz\"ahlbar unendlichen, symmetrischen Gruppe.}
Math. Z. 85, 40--61 (1964).

\bibitem{VK}
Vershik, A. M.; Kerov, S. V.
{\it Characters and factor representations
 of the infinite symmetric group.}
Soviet Math. Dokl. 23 (1981), no. 2, 389--392.

\bibitem{VK2}
Vershik, A. M.; Kerov, S. V.
{\it 
 Asymptotic theory of the characters of a symmetric group.}
Functional Anal. Appl. 15 (1981), no. 4, 246--255 (1982).

\bibitem{Wal}
Walsh, T.R.S. {\it Hypermaps versus bipartite maps,}
 J. Combinat. Theory B,
1975, vol. 18, 155--163


\end{thebibliography}
\end{document}